\documentstyle[12pt]{amsart}

\newsymbol\boxtimes 1202

\newcommand{\nc}{\newcommand}

\nc{\ad}{{\mbox{\bf{ad}}}}
\nc{\AJ}{{\operatorname{aj}}}
\nc{\Aut}{{\operatorname{Aut}}}
\nc{\Bls}{{{\cal B}ls}}
\nc{\Boxtimes}{{\fbox{$\times$}}} 
\nc{\blt}{{\bullet}}
\nc{\bSt}{{\mbox{\bf{St}}}}
\nc{\card}{{\operatorname{card}}}
\nc{\Cch}{{\check{C}}}
\nc{\cd}{{\operatorname{cd}}}
\nc{\Ch}{{\operatorname{Ch}}}
\nc{\chara}{{\operatorname{char}}}
\nc{\CHom}{{\cal{H}om}}
\nc{\Coker}{{\operatorname{Coker}}}
\nc{\codim}{{\operatorname{codim}}}
\nc{\Cone}{{\operatorname{Cone}}}
\nc{\cSgn}{{\cal{S}gn}}
\nc{\depth}{{\operatorname{depth}}} 
\nc{\dirlim}{{\underset{\rightarrow}{\operatorname{lim}}}}
\nc{\dotbox}{{\overset{\bullet}{\boxtimes}}}
\nc{\dotimes}{{\overset{\bullet}{\otimes}}}
\nc{\Ed}{{\operatorname{Edge}}}
\nc{\emp}{{\emptyset}}
\nc{\Ext}{{\operatorname{Ext}}}
\nc{\Fac}{{\cal{F}ac}}
\nc{\Fun}{{\operatorname{F}}}
\nc{\FS}{{\cal{FS}}}
\nc{\Hom}{{\operatorname{Hom}}}
\nc{\had}{{{\hat{\mbox{\bf{ad}}}}}}
\nc{\hgt}{{\operatorname{ht}}}
\nc{\Id}{{\operatorname{Id}}}
\nc{\id}{{\operatorname{id}}}
\nc{\Ima}{{\operatorname{Im}}}
\nc{\ind}{{\operatorname{ind}}}
\nc{\Ind}{{\operatorname{Ind}}}
\nc{\infi}{{\operatorname{inf}}}
\nc{\infh}{{\frac{\infty}{2}}}
\nc{\invlim}{{\underset{\leftarrow}{\operatorname{lim}}}}
\nc{\Jac}{{{\cal J}ac}}
\nc{\Ker}{{\operatorname{Ker}}}
\nc{\lcm}{{\operatorname{lcm}}}
\nc{\Locsys}{{{\cal L}ocsys}}
\nc{\Map}{{{\cal M}ap}}
\nc{\modul}{{\operatorname{mod}}}
\nc{\Mor}{{\operatorname{Mor}}}
\nc{\MS}{{\cal{MS}}}
\nc{\Ob}{{\operatorname{Ob}}}
\nc{\opp}{{\operatorname{opp}}}
\nc{\Or}{{{\cal O}r}}
\nc{\Ord}{{{\cal O}rd}}
\nc{\Part}{{{\cal P}art}}
\nc{\PGL}{{\operatorname{PGL}}}
\nc{\Pic}{{\operatorname{Pic}}}
\nc{\Rep}{{{\cal{R}}ep}}
\nc{\rk}{{\operatorname{rk}}}
\nc{\Sets}{{{\cal{S}}ets}}
\nc{\Sew}{{{\cal{S}}ew}}
\nc{\sgn}{{\operatorname{sgn}}}
\nc{\Sh}{{{\cal S}h}}
\nc{\Sign}{{{\cal S}ign}}
\nc{\Spe}{{\mbox{\bf{Sp}}}}
\nc{\supr}{{\operatorname{sup}}}
\nc{\Supp}{{\operatorname{Supp}}}
\nc{\supp}{{\operatorname{supp}}}
\nc{\Teich}{{{\cal{T}}eich}}
\nc{\tFS}{{\widetilde{\cal{FS}}}}
\nc{\Tor}{{\operatorname{Tor}}}
\nc{\totimes}{{\tilde{\otimes}}}
\nc{\tr}{{\operatorname{tr}}}
\nc{\tRep}{{\widetilde{{\cal R}ep}}}
\nc{\tTeich}{{\widetilde{{\cal T}eich}}}
\nc{\Vect}{{{\cal V}ect}}
\nc{\Ve}{{\operatorname{Vert}}}
\nc{\wt}{{\widetilde}}

\nc{\bo}{{\mbox{\bf{0}}}}
\nc{\One}{{\mbox{\bf{1}}}}
\nc{\one}{{\mbox{\bf{1}}}}
  
\nc{\BA}{{\Bbb A}}
\nc{\bA}{{\overline{A}}}
\nc{\ba}{{\mbox{\bf{a}}}}
\nc{\baB}{{\overline{B}}}
\nc{\baeta}{{\bar{\eta}}}
\nc{\baJ}{{\bar{J}}}
\nc{\BB}{{\Bbb B}}
\nc{\bB}{{\mbox{\bf{B}}}}
\nc{\bc}{{\mbox{\bf{c}}}}
\nc{\bC}{{\overline{C}}}
\nc{\BC}{{\Bbb{C}}}
\nc{\bCC}{{\overline{\cal{C}}}}
\nc{\bCM}{{\overline{\cal{M}}}}
\nc{\bD}{{\bar{D}}}
\nc{\BD}{{\overline{D}}}
\nc{\bd}{{\mbox{\bf{d}}}}
\nc{\BE}{{\overline{E}}}
\nc{\BF}{{\overline{F}}}
\nc{\bF}{{\mbox{\bf{F}}}}
\nc{\bg}{{\mbox{\bf{g}}}}
\nc{\bGamma}{{\overline{\Gamma}}}
\nc{\bL}{{\mbox{\bf{L}}}}
\nc{\blambda}{{\bar{\lambda}}}
\nc{\bM}{{\mbox{\bf{M}}}}
\nc{\bmu}{{\vec{\mu}}}
\nc{\BN}{{\Bbb{N}}}
\nc{\bnu}{{\mbox{\boldmath{${\nu}$}}}}
\nc{\bof}{{\mbox{\bf{f}}}}
\nc{\BP}{{\Bbb P}}
\nc{\bP}{{\mbox{\bf{P}}}}
\nc{\BPO}{{\overset{\circ}{\BP}}}
\nc{\BQ}{{\Bbb Q}}
\nc{\bq}{{\mbox{\bf{q}}}}
\nc{\BR}{{\Bbb{R}}}
\nc{\bR}{{\mbox{\bf{R}}}}
\nc{\br}{{\mbox{\bf{r}}}}
\nc{\breta}{{\bar{\eta}}}
\nc{\bs}{{\mbox{\bf{s}}}}
\nc{\bS}{{\mbox{\bf{S}}}}
\nc{\bt}{{\mbox{\bf{t}}}}
\nc{\bU}{{\mbox{\bf{U}}}}
\nc{\bu}{{\mbox{\bf{u}}}} 
\nc{\BUpsilon}{{\bar{\Upsilon}}}
\nc{\bw}{{\mbox{\bf{w}}}}
\nc{\bx}{{\mbox{\bf{x}}}}
\nc{\BZ}{{\Bbb{Z}}}
\nc{\bz}{{\mbox{\bf{z}}}}
\nc{\bzero}{\mbox{\boldmath{$0$}}}

\nc{\CA}{{\cal A}}
\nc{\CAD}{{\overset{\bullet}{\cal{A}}}} 
\nc{\CAO}{{\overset{\circ}{\cal{A}}}}
\nc{\CB}{{\cal B}}
\nc{\CC}{{\cal C}}
\nc{\CD}{{\cal D}}
\nc{\CE}{{\cal E}}
\nc{\CF}{{\cal F}}
\nc{\CH}{{\cal H}}
\nc{\CI}{{\cal I}}
\nc{\CID}{{\overset{\bullet}{\cal{I}}}}
\nc{\CJ}{{\cal J}}
\nc{\CK}{{\cal K}}
\nc{\CL}{{\cal L}}
\nc{\CM}{{\cal M}}
\nc{\CN}{{\cal N}}
\nc{\CO}{{\cal O}}
\nc{\CP}{{\cal P}}
\nc{\CPO}{{\overset{\circ}{\cal{P}}}}
\nc{\CQ}{{\cal Q}}
\nc{\CR}{{\cal R}}
\nc{\CS}{{\cal S}}
\nc{\CT}{{\cal T}}
\nc{\CTD}{{\overset{\bullet}{\cal{T}}}}
\nc{\CTPO}{{\overset{\circ}{\cal{T}\cal{P}}}}
\nc{\CU}{{\cal{U}}}
\nc{\CV}{{\cal V}}
\nc{\CX}{{\cal X}}
\nc{\CY}{{\cal Y}}
\nc{\CZ}{{\cal Z}}

\nc{\dCL}{{\overset{\bullet}{\cal{L}}}}
\nc{\dd}{{\operatorname{d}}}
\nc{\ddelta}{{\overset{\bullet}{\delta}}}
\nc{\dfu}{{\overset{\bullet}{\frak{u}}}}
\nc{\dlambda}{{\overset{\bullet}{\lambda}}}
\nc{\DO}{{\overset{\circ}{D}}}
\nc{\dpar}{{\partial}}
\nc{\dS}{{\overset{\bullet}{S}}}
\nc{\dT}{{\overset{\bullet}{T}}}

\nc{\fA}{{\frak{A}}}
\nc{\fC}{{\frak{C}}}
\nc{\fD}{{\frak{D}}}
\nc{\fE}{{\frak{E}}}
\nc{\fF}{{\frak{F}}}
\nc{\ff}{{\frak{f}}}
\nc{\fg}{{\frak{g}}}
\nc{\fH}{{\frak{H}}}
\nc{\fl}{{\frak{l}}}
\nc{\fL}{\frak{L}}
\nc{\fM}{\frak{M}}
\nc{\fp}{{\frak{p}}}
\nc{\fu}{{\frak{u}}}

\nc{\hCH}{{\hat{\cal{H}}}}
\nc{\hCI}{{\hat{\cal{I}}}}
\nc{\hfC}{{\hat{\frak{C}}}}
\nc{\hfg}{{\hat{\frak{g}}}}
\nc{\hL}{{\hat{L}}}
\nc{\HO}{{\overset{\circ}{H}}}
\nc{\hpsi}{{\hat{\psi}}}
\nc{\hx}{{\hat{x}}}

\nc{\jo}{{\overset{\circ}{j}}}

\nc{\phid}{{\overset{\bullet}{\phi}}}

\nc{\tA}{{\tilde{A}}}
\nc{\ta}{{\tilde{a}}}
\nc{\tB}{{\tilde{B}}}
\nc{\tb}{{\tilde{b}}} 
\nc{\tBP}{{\tilde{\BP}}}
\nc{\tC}{{\tilde{C}}}
\nc{\tc}{{\tilde{c}}}
\nc{\tCA}{{\tilde{\cal{A}}}}
\nc{\tCC}{{\tilde{\cal{C}}}}
\nc{\tCH}{{\tilde{\cal{H}}}}
\nc{\tCI}{{\tilde{\cal{I}}}}
\nc{\tCO}{{\tilde{\cal{O}}}}
\nc{\tCP}{{\tilde{\cal{P}}}}
\nc{\tCT}{{\tilde{\cal{T}}}}
\nc{\tD}{{\tilde{D}}}
\nc{\tDelta}{{\tilde{\Delta}}}
\nc{\tE}{{\tilde E}}
\nc{\tF}{{\tilde F}}
\nc{\tfD}{{\tilde{\frak{D}}}}
\nc{\tfF}{{\tilde{\frak{F}}}}
\nc{\tff}{{\tilde{\frak{f}}}}
\nc{\tfu}{{\tilde{\frak{u}}}}
\nc{\tJ}{{\tilde{J}}}
\nc{\tj}{{\tilde{j}}}
\nc{\tK}{{\tilde K}}
\nc{\tL}{{\tilde{L}}}
\nc{\tM}{{\tilde{M}}}
\nc{\tP}{{\tilde{P}}}
\nc{\tPhi}{{\tilde{\Phi}}}
\nc{\tpi}{\tilde{\pi}}
\nc{\TPO}{{\overset{\circ}{T\BP}}}
\nc{\tR}{{\tilde{R}}}
\nc{\tS}{{\tilde S}}
\nc{\tT}{{\tilde{T}}}
\nc{\ttau}{{\tilde{\tau}}}
\nc{\ttheta}{{\tilde{\theta}}}
\nc{\tU}{{\tilde{U}}}
\nc{\tUpsilon}{{\tilde{\Upsilon}}}
\nc{\ty}{{\tilde y}}
\nc{\tY}{{\tilde Y}}
\nc{\txi}{{\tilde{\xi}}}

\nc{\UD}{{\overset{\bullet}{U}}}
\nc{\UO}{{\overset{\circ}{U}}}

\nc{\vA}{{\vec{A}}}
\nc{\valpha}{{\vec{\alpha}}}
\nc{\vbeta}{{\vec{\beta}}}
\nc{\vc}{{\vec{c}}}
\nc{\vD}{{\vec{D}}}
\nc{\vd}{{\vec{d}}}
\nc{\vgamma}{{\vec{\gamma}}}
\nc{\vK}{{\vec{K}}}
\nc{\vlambda}{{\vec{\lambda}}}
\nc{\vmu}{{\vec{\mu}}}
\nc{\vnu}{{\vec{\nu}}}
\nc{\vo}{{\vec{0}}}
\nc{\vu}{{\vec{u}}}
\nc{\vx}{{\vec{x}}}
\nc{\vy}{\vec{y}}
\nc{\vzero}{\vec{0}}

\nc{\XO}{{\overset{\circ}{X}}}

\nc{\ya}{{\operatorname{aj}}}

\nc{\nen}{\newenvironment}
\nc{\ol}{\overline}
\nc{\ul}{\underline}
\nc{\ra}{\rightarrow}
\nc{\lra}{\longrightarrow}
\nc{\Lra}{\Longrightarrow}
\nc{\lla}{\longleftarrow}
\nc{\Llra}{\Longleftrightarrow}
\nc{\hra}{\hookrightarrow}
\nc{\iso}{\overset{\sim}{\lra}}
\nc{\rlh}{\rightleftharpoons} 

\nc{\Thm}[1]{Theorem~\ref{#1}}
\nc{\Prop}[1]{Proposition~\ref{#1}}
\nc{\Lem}[1]{Lemma~\ref{#1}}
\nc{\Cor}[1]{Corollary~\ref{#1}}
\nc{\Conj}[1]{Conjecture~\ref{#1}}
\nc{\Claim}[1]{Claim~\ref{#1}}
\nc{\Defn}[1]{Definition~\ref{#1}}
\nc{\Exa}[1]{Example~\ref{#1}}
\nc{\Rem}[1]{Remark~\ref{#1}}
\nc{\Note}[1]{Note~\ref{#1}}

\nen{thm}[1]{\label{#1}{\bf Theorem.\ } \em}{}
\nen{prop}[1]{\label{#1}{\bf Proposition.\ } \em}{}
\nen{lem}[1]{\label{#1}{\bf Lemma.\ } \em}{}
\nen{cor}[1]{\label{#1}{\bf Corollary.\ } \em}{}
\nen{conj}[1]{\label{#1}{\bf Conjecture.\ } \em}{}
\nen{claim}[1]{\label{#1}{\bf Claim.\ } \em}{}

\nen{defn}[1]{\label{#1}{\bf Definition.\ } }{}
\nen{exa}[1]{\label{#1}{\bf Example.\ } }{}

\nen{rem}[1]{\label{#1}{\em Remark.\ } }{}
\nen{note}[1]{\label{#1}{\em Note.\ } }{}
\nen{exer}[1]{\label{#1}{\em Exercise.\ } }{}

\begin{document}

\title[]{Localization of modules over small quantum groups}
\author{Michael Finkelberg}
\address{M.F.: Independent Moscow University, Miklukho-Maklai St., 65-3-86,  
Moscow 117342 Russia} 
\email{fnklberg@@ium.ips.ras.ru}
\author{Vadim Schechtman}
\address{V.S.: 26 Bakinskikh Komissarov St., 3-1-22, Moscow 117571 Russia}
\thanks{The first author was partially supported by
grants from AMS, SOROS, INTAS. 
The second author was partially supported by NSF.}
\date{March 1996\\
q-alg/9604001}
\maketitle


\section{Introduction}

This article is a report on \cite{fs}, \cite{bfs}.   

\subsection{} Let $(I,\cdot)$ be a Cartan datum of finite type and  
$(Y,X,\ldots)$ the simply connected root datum of type $(I,\cdot)$, 
cf. ~\cite{l1}.  

Let $l>1$ be an integer. Set $\ell=l/(l,2)$; to simplify the exposition,  
we will suppose in this Introduction that $d_i:=i\cdot i/2$ divides $\ell$ 
for all $i\in I$; we set $\ell_i:=\ell/d_i$. We suppose that all 
$\ell_i>3$.  
Let $\rho\in X$ be the half-sum 
of positive roots; let $\rho_{\ell}\in X$ be defined by $\langle i,\rho_{\ell}
\rangle =\ell_i-1\ (i\in I)$. We define a lattice $Y_{\ell}\subset Y$ by 
$Y_{\ell}=\{\mu\in X|\ \mbox{for all}\ \mu'\in X, \mu\cdot\mu'\in\ell\BZ\}$, 
and set $\dd_{\ell}=\card(X/Y_{\ell})$.    

We fix a base field $k$ of characteristic not dividing $l$. We suppose that 
$k$ contains a primitive $l$-th root of unity $\zeta$, and fix it. 
Starting from these data, one defines certain category 
$\CC$. Its objects are finite dimensional $X$-graded $k$-vector spaces 
equipped with an action of Lusztig's "small" quantum group $\fu$ 
(cf. ~\cite{l2}) such that the action of its Cartan subalgebra is compatible 
with the $X$-grading. Variant: one defines certain algebra $\dfu$ which is an  
"$X$-graded" version of $\fu$ (see \ref{dot}), and an object of 
$\CC$ is a finite 
dimensional $\dfu$-module.  For the precise definition of $\CC$, see  
~\ref{def c}, \ref{comp lu}. 
For $l$ prime and $\chara k=0$, the category $\CC$ was studied in ~\cite{ajs}, 
for $\chara k>0$ and arbitrary $l$, $\CC$ was studied in \cite{aw}.  
The category $\CC$ 
admits a remarkable structure of a {\em ribbon}\footnote{in other terminology, 
braided balanced rigid tensor} category (Lusztig). 

\subsection{} 
The main aim of this work is to introduce      
certain tensor category $\FS$ of geometric origin, which is  
equivalent to $\CC$. Objects of $\FS$ are called {\bf (finite) factorizable 
sheaves}. 

A notion of a factorizable sheaf is the main new concept of this paper. 
Let us give an informal idea, what kind of an object is it. 
Let $D$ be the unit disk in the complex plane\footnote{One could also take   
a complex affine line or a formal disk; after a suitable modification 
of the definitions, the resulting categories of factorizable sheaves 
are canonically equivalent.}. Let 
$\CD$ denote the space of positive $Y$-valued divisors on $D$. Its points 
are formal linear combinations $\sum\nu\cdot x\ (\nu\in Y^+:=\BN[I], x\in D)$. 
This space is a disjoint union 
$$
\CD=\coprod_{\nu\in Y^+}\ D^{\nu}
$$
where $D^{\nu}$ is the subspace of divisors of degree $\nu$. Variant: $D^{\nu}$ 
is the configuration space of $\nu$ points running on $D$; its points are 
(isomorphism clases of) pairs of maps 

($J\lra D$, $j\mapsto x_j$;\ $\pi: J\lra I$, $j\mapsto i_j$),
 
$J$ being a finite set.  
We say that we have a finite set $\{x_j\}$ of (possibly equal) 
points of $D$, a point $x_j$ being "coloured" by the colour $i_j$. The sum  
(in $Y$) of colours of all points should be equal to $\nu$. The space 
$D^{\nu}$ is a smooth affine analytic variety; it carries a canonical 
stratification defined by various intersections of hypersurfaces $\{x_j=0\}$ and 
$\{x_{j'}=x_{j''}\}$. The open stratum (the complement of all the hypersurfaces 
above) is denoted by $A^{\nu\circ}$.     

One can imagine $\CD$ as a $Y^+$-graded smooth stratified variety. 
Let $\CA=\coprod A^{\nu}\subset\CD$ be the open $Y^+$-graded subvariety  
of positive $Y$-valued divisors on $D-\{0\}$. We have an open $Y^+$-graded 
subvariety $\CA^{\circ}=\coprod A^{\nu\circ}\subset\CA$. 

Let us consider the $(Y^+)^2$-graded variety  
$$
\CA\times\CA=\coprod\ A^{\nu_1}\times A^{\nu_2}; 
$$
we define another $(Y^+)^2$-graded variety $\widetilde{\CA\times\CA}$, 
together with two maps 

(a) $\CA\times\CA\overset{p}{\lla}\widetilde
{\CA\times\CA}\overset{m}{\lra}\CA$  

respecting the $Y^+$-gradings\footnote{a $(Y^+)^2$-graded space is considered 
as a $Y^+$-graded by means of the addition $(Y^+)^2\lra Y^+$.}; the map 
$p$ is a homotopy equivalence. One can imagine the   
diagram (a) above as a "homotopy multiplication"    
$$
m_{\CA}:\CA\times\CA\lra\CA; 
$$
this "homotopy map" is "homotopy associative"; the meaning of this is explained 
in the main text. We say that 
$\CA$ is a ($Y^+$-graded) "homotopy monoid"; $\CA^{\circ}\subset\CA$ is a 
"homotopy submonoid".  

The space $\CD$ is a "homotopy $\CA$-space": there is a    
"homotopy map"   
$$
m_{\CD}:\CD\times\CA\lra\CD
$$
which is, as above, a diagram of usual maps between $Y^+$-graded 
varieties    

(b) $\CD\times \CA\overset{p}{\lla}\widetilde{\CD\times\CA}
\overset{m}{\lra}\CD$,   

$p$ being a homotopy equivalence.  

For each $\mu\in X,\ \nu\in Y^+$, one defines a one-dimensional 
$k$-local system 
$\CI_{\mu}^{\nu}$ over $A^{\nu\circ}$. Its monodromies are defined by variuos 
scalar products of "colours" of running points and the colour $\mu$ of the 
origin $0\in D$.  
These local systems have the following 
compatibility. For each $\mu\in X, \nu_1,\nu_2\in Y^+$, a 
{\em factorization isomorphism}  
$$
\phi_{\mu}(\nu_1,\nu_2): m^*\CI_{\mu}^{\nu_1+\nu_2}\iso p^*(\CI^{\nu_1}_{\mu}
\boxtimes\CI^{\nu_2}_{\mu-\nu_1})
$$
is given (where $m, p$ are the maps in the diagram (a) above), 
these isomorphisms satisfying certain {\em (co)associativity property}. The 
collection of local systems $\{\CI^{\nu}_{\mu}\}$ is an $X$-graded 
local system $\CI=\oplus\ \CI_{\mu}$ over $\CA^{\circ}$. One could imagine 
the collection of factorization isomorphisms $\{\phi_{\mu}(\nu_1,\nu_2)\}$ 
as an isomorphism 
$$
\phi:\ m_{\CA^{\circ}}^*\CI\iso\CI\boxtimes\CI
$$
We call $\CI$ the {\em braiding local system}. 

Let $\CI^{\bullet}$ be the perverse sheaf on $\CA$ which is the 
Goresky-MacPherson 
extension of the local system $\CI$. The isomorphism $\phi$ above induces 
a similar isomorphism 
$$
\phi^{\bullet}:\ m_{\CA}^*\CI^{\bullet}\iso\CI^{\bullet}\boxtimes\CI^{\bullet} 
$$      
This sheaf, together with $\phi^{\bullet}$, looks like a "coalgebra"; it  
is an incarnation of the quantum group $\dfu$.  

A factorizable sheaf is a couple 

($\mu\in X$ ("the highest weight"); 
a perverse sheaf $\CX$ over $\CD$, smooth along the canonical stratification). 
 
Thus, $\CX$ is a collection of sheaves $\CX^{\nu}$ over $\CD^{\nu}$. These 
sheaves should be connected by {\em factorization isomorphisms}
$$
\psi(\nu_1,\nu_2):\ m^*\CX^{\nu_1+\nu_2}\iso p^*(\CX^{\nu_1}
\boxtimes\CI^{\nu_2\bullet}_{\mu-\nu_1})
$$
satisfying an associativity property. Here $m, p$ are as in the diagram (b) 
above. One could imagine the whole collection $\{\psi(\nu_1,\nu_2)\}$ as  
an isomorphism 
$$
\psi: m_{\CD}^*\CX\iso\CI^{\bullet}\boxtimes\CX
$$
satisfying a (co)associativity property. We impose also certain finiteness 
(of singularities) condition on $\CX$. So, this object looks like a 
"comodule" over $\CI^{\bullet}$. It is an incarnation of an $\dfu$-module. 

We should mention one more important part of the structure of the space 
$\CD$. It comes from natural closed embeddings 
$\iota_{\nu}:\ \CD\hra\CD[\nu]$ (were $[\nu]$ denotes the 
shift of the grading); 
these mappings define certain inductive system, and  
a factorizable sheaf is a sheaf on its inductive limit. 

The latter inductive 
limit is an example of a {\bf "semiinfinite space"}.     

For the precise definitions, see Sections \ref{braiding}, \ref{factoriz} 
below. 
  
The category $\FS$ has a structure of a  
braided balanced tensor category coming from geometry. The  
tensor structure on $\FS$ is defined using the functors of nearby cycles. 
The tensor equivalence 
$$
\Phi:\FS\iso\CC
$$
is defined using vanishing cycles functors. It respects the  
braidings and balances. This is the contents of Part I.  

\subsection{} 
\label{glue intro} Factorizable sheaves are local objects. It turns out that    
one can "glue" them along complex curves. More precisely,   
given a finite family $\{\CX_a\}_{a\in A}$ of factorizable sheaves 
and a smooth proper curve $C$ together 
with a family $\{x_a,\tau_a\}_A$ of distinct points $x_a\in C$ with 
non-zero tangent vectors $\tau_a$ at them (or, more generally, a family of such 
objects over a smooth base $S$), one can define a perverse 
sheaf, denoted by   
$$
\Boxtimes_A^{(C)}\ \CX_a
$$
on the (relative) configuration space $C^{\nu}$. Here $\nu\in Y^+$ is 
defined by 

(a) $\nu=\sum_{a\in A}\ \mu_a+(2g-2)\rho_{\ell}$ 

where $g$ is the genus of $C$, $\mu_a$ is the highest weight of $\CX_a$. 
We {\em assume} that the right hand side of the equality belongs to $Y^+$.   

One can imagine this sheaf as an "exterior tensor product" 
of the family $\{\CX_a\}$ along $C$. It is obtained by "planting" the 
sheaves $\CX_a$ into the points $x_a$. To glue them together, one needs 
a "glue". 
This glue is called the {\bf Heisenberg local system} $\CH$; it was 
constructed by R.Bezrukavnikov. $\CH$ is a sister 
of the braiding local system $\CI$. Let us describe what it is. 

For a finite set $A$, let $\CM_A$ denote the moduli stack of 
punctured curves $(C,\{x_a,\tau_a\}_A)$ as above; let $\eta: C_A\lra\CM_A$ be 
the universal curve. Let $\lambda_A=\det (R\eta_*\CO_{C_A})$ be the determinant 
line bundle on $\CM_A$, and $\CM_{A;\lambda}\lra\CM_A$ be the total space 
of $\lambda_A$ with the zero section removed. For $\nu\in Y^+$, let 
$\eta^{\nu}: C^{\nu}_A\lra\CM_A$ be the relative configuration space 
of $\nu$ points running on $C_A$; let 
$C^{\nu\circ}_A\subset C^{\nu}_A$ be the 
open stratum (where the running points are distinct and distinct from 
the punctures $x_a$). 
The complementary subscript $()_{\lambda}$ will denote the base 
change under $\CM_{A;\lambda}\lra\CM_A$. The complementary subscript $()_g$ 
will denote the base change under $\CM_{A,g}\hra\CM_A$, $\CM_{A,g}$ being 
the substack of curves of genus $g$. 

The Heisenberg local system is a collection of local systems 
$\CH_{\vmu;A,g}$ over the stacks $C^{\nu\circ}_{A,g;\lambda}$. Here 
$\vmu$ is an $A$-tuple $\{\mu_a\}\in X^a$; $\nu=\nu(\vmu;g)$ is defined 
by the equality (a) above. We assume that $\vmu$ is such that the right 
hand side of this equality really belongs to $Y^+$. The dimension of 
$\CH_{\vmu;A,g}$ is equal to $\dd_{\ell}^g$; the monodromy around the 
zero section of the determinant line bundle is equal to 

(b) $c=(-1)^{\card(I)}\zeta^{-12\rho\cdot\rho}$.

This formula is due to R.Bezrukavnikov. 
These local systems have a remarkable compatibility ("fusion") property 
which we do not specify here, see Section \ref{fusion}.        

\subsection{}    
In Part II we study the sheaves $\Boxtimes_A^{(C)}\ \CX_a$ for  
$C=\BP^1$. Their cohomology, when expressed algebraically,  
turn out to coincide with certain "semiinfinite" $\Tor$ spaces in  
the category $\CC$ introduced by S.M.Arkhipov. Due to the results  
of Arkhipov, this enables one to prove that  
"spaces of conformal blocks" in WZW models are the natural subquotients of 
such cohomology spaces. 

\subsection{} 
In Part III we study the sheaves $\Boxtimes_A^{(C)}\ \CX_a$ for arbitrary 
smooth families of punctured curves. Let $\Boxtimes_{A,g}\ \CX_a$ denotes 
the "universal exterior product"  
living on $C^{\nu}_{A,g;\lambda}$ where $\nu=\nu(\vmu;g)$ is as 
in \ref{glue intro} (a) above, 
$\vmu=\{\mu_a\},\ \mu_a$ being the highest weight of $\CX_a$. Let us integrate 
it: consider 
$$
\int_{C^{\nu}}\ \Boxtimes_{A,g}\ \CX_a:=R\eta^{\nu}_*(\Boxtimes_{A,g}\ \CX_a);
$$
it is a complex of sheaves over $\CM_{A,g;\lambda}$ with smooth cohomology; 
let us denote it $\langle\otimes_A\ \CX_a\rangle_g$. The cohomology sheaves 
of such complexes     
define a {\em fusion structure}\footnote{actually 
a family of such structures (depending on the degree of cohomology)} on $\FS$ 
(and hence on $\CC$). 
 
The classical WZW model fusion category 
is, in a certain sense, a subquotient of one of them. The number $c$, 
\ref{glue intro} (b),  
coincides with the "multiplicative central charge" of the model (in the 
sense of \cite{bfm}).  
 
As a consequence of this geometric description and of the Purity Theorem, 
\cite{bbd}, 
the local systems of conformal blocks (in arbitrary genus)  
are semisimple. The Verdier duality induces a canonical non-degenerate 
Hermitian form on them (if $k=\BC$). 

Almost all the results of Part III are due to R.Bezrukavnikov, cf. \cite{bfs}. 

\subsection{} The works \cite{fs} and \cite{bfs} in a sense complete 
the program outlined in \cite{s} (no fusion structure in higher genus 
was mentioned there).  We should mention a very interesting 
related work of G.Felder and collaborators.  
The idea of realizing quantum 
groups' modules in the cohomology of configuration spaces appeared 
independently in \cite{fw}. The Part III of the present work  
arose from our attempts to understand \cite{cfw}. 

We are grateful to D.Kazhdan for the encouragement. 
The first-named author was strongly influenced by the
ideas and suggestions of A.Beilinson, P.Deligne, B.Feigin, G.Felder and 
D.Kazhdan; he feels 
greatly indebted to all of them.  We are very obliged to 
G.Lusztig, and especially to S.M.Arkhipov and 
R.Bezrukavnikov, for the permission to use their unpublished
results.    

\subsection{} This article contains almost all main definitions and theorems  
from \cite{fs}, \cite{bfs}, but no proofs; we refer the  
reader to {\em op. cit.} for them. However, our exposition here 
differs from {\em op. cit.}; we hope it clarifies the subject to some extent.   
  

\newpage
\begin{center}
{\bf Part I. Local.}
\end{center} 

\section{The category $\CC$} 

\subsection{} 
\label{root} We will follow Lusztig's terminology 
and notations concerning root systems, cf. ~\cite{l1}.  

We fix an irreducible Cartan datum $(I,\cdot)$ 
of finite type. Thus, $I$ is a finite set together with a symmetric 
$\BZ$-valued bilinear form $\nu_1,\nu_2\mapsto\nu_1\cdot\nu_2$ 
on the free abelian group $\BZ[I]$, such that for any $i\in I$, 
$i\cdot i$ is even and positive, for any $i\neq j$ in $I$, 
$2i\cdot j/i\cdot i\leq 0$ and the $I\times I$-matrix $(i\cdot j)$ is 
positive definite. We set $d_i=i\cdot i/2$.

Let $d=\max_{i\in I}\ d_i$. This number is equal to the least 
common multiple of the numbers $d_i$, and belongs to the set $\{1,2,3\}$. 
For a simply laced root datum $d=1$.  

We set $Y=\BZ[I], X=\Hom(Y,\BZ)$; $\langle,\rangle:\ Y\times X\lra 
\BZ$ will denote the obvious pairing. The obvious 
embedding $I\hra Y$ will be denoted by $i\mapsto i$. We will denote by 
$i\mapsto i'$  the embedding 
$I\hra X$ given by $\langle i,j'\rangle=2i\cdot j/i\cdot i$. 
(Thus, $(Y,X,\ldots)$ is the simply connected root datum of type $(I,\cdot)$, 
in the terminology of Lusztig.) 

The above embedding $I\subset X$ extends by additivity to the embedding 
$Y\subset X$. We will regard $Y$ as the sublattice of $X$ by means 
of this embedding.  
For $\nu\in Y$, we will denote by the same letter $\nu$ 
its image in $X$. We set $Y^+=\BN[I]\subset Y$. 

We will use the following partial order on $X$. For $\mu_1, \mu_2\in X$, 
we write $\mu_1\leq \mu_2$ if $\mu_2-\mu_1$ belongs to $Y^+$.    

\subsection{} 
\label{l} We fix a base field $k$, an integer $l>1$ and a primitive root 
of unity $\zeta\in k$ as in the Introduction.
 
We set $\ell=l$ if $l$ is odd and $\ell=l/2$ 
if $l$ is even. For $i\in I$, we set $\ell_i=\ell/(\ell,d_i)$. Here 
$(a,b)$ stands for the greatest common divisor of $a,b$. We set 
$\zeta_i=\zeta^{d_i}$. We will assume that $\ell_i>1$ for any $i\in I$ and 
$\ell_i>-\langle i,j'\rangle+1$ for any $i\neq j$ in $I$. 

We denote by $\rho$ (resp. $\rho_{\ell}$) the element of $X$ such that 
$\langle i,\rho\rangle=1$ (resp.  
$\langle i,\rho_{\ell}\rangle=\ell_i-1$) for all $i\in I$. 

For a coroot $\beta\in Y$, there exists an element $w$ of the Weyl 
group $W$ of our Cartan datum  
and $i\in I$ such that $w(i)=\beta$. We set $\ell_{\beta}:=\frac{\ell}
{(\ell,d_i)}$; this number does not depend on the choice of $w$ and $i$. 
 
We have $\rho=\frac{1}{2}\sum\ \alpha;\ \rho_{\ell}=\frac{1}{2}\sum\ 
(\ell_{\alpha}-1)\alpha$, the sums over all positive roots $\alpha\in X$.  

For $a\in\BZ, i\in I$, we set $[a]_i=1-\zeta_i^{-2a}$.

\subsection{} 
\label{lattice} We use the same notation $\mu_1,\mu_2\mapsto\mu_1\cdot\mu_2$ 
for the unique extension of the bilinear form on $Y$ to a $\Bbb{Q}$-valued 
bilinear form on $Y\otimes_{\BZ}\Bbb{Q}=X\otimes_{\BZ}\Bbb{Q}$. 

We define a lattice $Y_{\ell}=\{\lambda\in X|\mbox{ for all }\mu\in X, 
\lambda\cdot\mu\in\ell\BZ\}$.  

\subsection{}
Unless specified otherwise, 
a "vector space" will mean a vector space over $k$; $\otimes$ will 
denote the tensor product over $k$. A "sheaf" (or a "local system")  
will mean a sheaf (resp. local system) of $k$-vector spaces.  

If $(T,\CS)$ is an open subspace of   
the space of complex points of a separate scheme of 
finite type over $\BC$, with the usual topology,  
together with an algebraic stratification $\CS$ satisfying the properties 
\cite{bbd} 2.1.13 b), c), we will denote by 
$\CM(T;\CS)$ the category of perverse sheaves over $T$ lisse along 
$\CS$, with respect to the middle perversity, cf. \cite{bbd} 2.1.13, 
2.1.16. 
     
\subsection{} Let $'\ff$ be the free associative $k$-algebra with $1$ 
with generators $\theta_i\ (i\in I)$. For $\nu=\sum \nu_ii\in\BN[I]$, 
let $'\ff_{\nu}$ be the subspace of $'\ff$ spanned by the monomials 
$\theta_{i_1} 
\cdot\ldots\cdot\theta_{i_a}$ such that $\sum_{j} i_j=\nu$ in $\BN[I]$.    

Let us regard $'\ff\otimes\ '\ff$ 
as a $k$-algebra with the product $(x_1\otimes x_2)(y_1\otimes y_2)= 
\zeta^{\nu\cdot\mu}x_1y_1\otimes x_2y_2\ (x_2\in\ '\ff_{\nu}, y_1\in\ 
'\ff_{\mu}$). Let $r$ denote a unique homomorphism of $k$-algebras 
$'\ff\lra\ '\ff\otimes\ '\ff$ carrying $\theta_i$ to $1\otimes\theta_i+
\theta_i\otimes 1\ (i\in I)$. 

\subsection{Lemma-definition} {\em There exists a unique $k$-valued bilinear 
form 
$(\cdot,\cdot)$ on $'\ff$ such that 

{\em (i)} $(1,1)=1;\ (\theta_i,\theta_j)=\delta_{ij}\ (i,j\in I)$; 
{\em (ii)} $(x,yy')=(r(x),y\otimes y')$ for all $x,y,y'\in\ '\ff$.  

This bilinear form is symmetric.} $\Box$  

In the right hand side of the equality (ii) we use the same notation 
$(\cdot,\cdot)$ for the bilinear form on $'\ff\otimes\ '\ff$ 
defined by $(x_1\otimes x_2,y_1\otimes y_2)=(x_1,y_1)(x_2,y_2)$. 

The radical of the form $'\ff$ is a two-sided ideal of $'\ff$. 

\subsection{} Let us consider the associative $k$-algebra $\fu$ 
(with $1$) defined by the generators $\epsilon_i, \theta_i\ (i\in I), 
K_{\nu}\ (\nu\in Y)$ and the relations (a) --- (e) below. 

(a) $K_0=1,\ K_{\nu}\cdot K_{\mu}=K_{\nu+\mu}\ (\nu,\mu\in Y)$; 

(b) $K_{\nu}\epsilon_i=\zeta^{\langle\nu,i'\rangle}\epsilon_iK_{\nu}\ 
(i\in I, \nu\in Y)$; 

(c) $K_{\nu}\theta_i=\zeta^{-\langle\nu,i'\rangle}\theta_iK_{\nu}\ 
(i\in I, \nu\in Y)$; 

(d) $\epsilon_i\theta_j-\zeta^{i\cdot j}\theta_j\epsilon_i=\delta_{ij}
(1-\tK_i^{-2})\ (i,j\in I)$. 

Here we use the notation $\tK_{\nu}=\prod_i\ K_{d_i\nu_ii}\ (\nu=\sum\nu_ii)$. 

(e) If $f(\theta_i)\in\ '\ff$ belongs to the radical of the form 
$(\cdot,\cdot)$ then $f(\theta_i)=f(\epsilon_i)=0$ in $\fu$. 

\subsection{} There is a unique $k$-algebra homomorphism $\Delta: \fu\lra 
\fu\otimes\fu$ such that 
$$
\Delta(\epsilon_i)=\epsilon_i\otimes 1+\tK_i^{-1}\otimes\epsilon_i;\ 
\Delta(\theta_i)=\theta_i\otimes 1+\tK_i^{-1}\otimes\theta_i;\ 
\Delta(K_{\nu})=K_{\nu}\otimes K_{\nu}
$$
for any $i\in I, \nu\in Y$. Here $\fu\otimes\fu$ is regarded as an algebra 
in the standard way. 

There is a unique $k$-algebra homomorphism $e:\fu\lra k$ 
such that $e(\epsilon_i)=e(\theta_i)=0, e(K_{\nu})=1\ (i\in I, \nu\in Y)$. 

\subsection{} There is a unique $k$-algebra homomorphism $A:\fu\lra\fu^{
\opp}$ such that 
$$
A(\epsilon_i)=-\epsilon_i\tK_i;\ A(\theta_i)=-\theta_i\tK_i;\ 
A(K_{\nu})=K_{-\nu}\ (i\in I, \nu\in Y).
$$
There is a unique $k$-algebra homomorphism $A':\fu\lra\fu^{\opp}$ such that 
$$
A'(\epsilon_i)=-\tK_i\epsilon_i;\ A'(\theta_i)=-\tK_i\theta_i;\ 
A'(K_{\nu})=K_{-\nu}\ (i\in I, \nu\in Y).
$$

\subsection{} The algebra $\fu$ together with the additional structure given 
by the comultiplication $\Delta$, the counit $e$, the antipode $A$ and 
the skew-antipode $A'$, is a Hopf algebra. 

\subsection{} 
\label{def c} Let us define a category $\CC$ as follows. An object 
of $\CC$ is a $\fu$-module $M$ which is finite dimensional over $k$, 
with a given direct sum decomposition $M=\oplus_{\lambda\in X}\ M_{\lambda}$ 
(as a vector space) such that $K_{\nu}x=\zeta^{\langle\nu,\lambda\rangle}
x$ for any $\nu\in Y, \lambda\in X, x\in M_{\lambda}$. A morphism in $\CC$ 
is a $\fu$-linear map respecting the $X$-gradings. 

Alternatively, an object of $\CC$ may be defined as an $X$-graded finite 
dimensional vector space $M=\oplus\ M_{\lambda}$ equipped with linear 
operators
$$
\theta_i: M_{\lambda}\lra M_{\lambda-i'},\ 
\epsilon_i: M_{\lambda}\lra M_{\lambda+i'}\ (i\in I, \lambda\in X)
$$
such that 

(a) for any $i,j\in I, \lambda\in X$, the operator 
$\epsilon_i\theta_j-\zeta^{i\cdot j}\theta_j\epsilon_i$ acts as the 
multiplication by 
$\delta_{ij}[\langle i,\lambda\rangle]_i$ on $M_{\lambda}$. 

Note that $[\langle i,\lambda\rangle]_i=
[\langle d_ii,\lambda\rangle]=[i'\cdot\lambda]$. 
 
(b) If $f(\theta_i)\in\ '\ff$ belongs to the radical of the form 
$(\cdot,\cdot)$ then the operators $f(\theta_i)$ and $f(\epsilon_i)$ act 
as zero on $M$. 

\subsection{} In ~\cite{l2} Lusztig defines an algebra $\bu_{\CB}$ over the 
ring $\CB$ which is a quotient of $\BZ[v,v^{-1}]$ ($v$ being an 
indeterminate) by the $l$-th cyclotomic polynomial. Let us consider 
the $k$-algebra $\bu_k$ obtained from $\bu_{\CB}$ by the base change 
$\CB\lra k$ sending $v$ to $\zeta$. The algebra $\bu_k$ is generated 
by certain elements $E_i, F_i, K_i\ (i\in I)$. Here $E_i=E_{\alpha_i}^{(1)},  
F_i=F_{\alpha_i}^{(1)}$ in the notations of {\em loc.cit.} 

Given an object $M\in\CC$, let us introduce the operators $E_i, F_i, K_i$ on 
it by 
$$
E_i=\frac{\zeta_i}{1-\zeta_i^{-2}}\epsilon_i\tK_i,\ F_i=\theta_i, K_i=K_i. 
$$

\subsection{Theorem} 
\label{comp lu} {\em The above formulas define the action 
of the Lusztig's algebra {\em $\bu_k$} on an object $M$. 

This rule defines an equivalence of $\CC$ with the category whose 
objects are $X$-graded finite dimensional {\em $\bu_k$}-modules 
$M=\oplus M_{\lambda}$ such that $K_ix=\zeta^{\langle i,\lambda\rangle}x$ 
for any $i\in I, \lambda\in X, x\in M_{\lambda}$.} $\Box$ 

\subsection{} The structure of a Hopf algebra on $\fu$ defines canonically 
a {\em rigid tensor} structure on $\CC$ (cf. \cite{kl}IV, Appendix).  

The Lusztig's algebra $\bu_k$ also has an additional structure of a Hopf 
algebra. It induces the same rigid tensor structure on $\CC$. 

We will denote the duality in $\CC$ by $M\mapsto M^*$. The unit 
object will be denoted by $\One$.  

\subsection{} Let $\fu^-$ (resp. $\fu^+$, $\fu^0$) denote the $k$-subalgebra 
generated by the elements $\theta_i\ (i\in I)$ (resp. $\epsilon_i\ 
(i\in I)$, $K_{\nu}\ (\nu\in Y)$). We have the triangular decomposition 
$\fu=\fu^-\fu^0\fu^+=\fu^+\fu^0\fu^-$. 

We define the "Borel" subalgebras $\fu^{\leq 0}=\fu^-\fu^0,\ \fu^{\geq 0}=
\fu^+\fu^0$; they are the Hopf subalgebras of $\fu$.   

Let us introduce the $X$-grading 
$\fu=\oplus\fu_{\lambda}$ as a unique grading compatible with 
the structure of an algebra such that $\theta_i\in\fu_{-i'},\  
\epsilon_i\in\fu_{i'},\ K_{\nu}\in\fu_0$. We will use the induced 
gradings on the subalgebras of $\fu$. 

\subsection{} 
\label{ind} Let $\CC^{\leq 0}$ (resp. $\CC^{\geq 0}$) be the category 
whose objects are $X$-graded finite dimensional $\fu^{\leq 0}$- 
(resp. $\fu^{\geq 0}$-) modules $M=\oplus M_{\lambda}$ such that 
$K_{\nu}x=\zeta^{\langle\nu,\lambda\rangle}x$ for any $\nu\in Y, 
\lambda\in X, x\in M_{\lambda}$. Morphisms are $\fu^{\leq 0}$- 
(resp. $\fu^{\geq 0}$-) linear maps compatible with the $X$-gradings. 

We have the obvious functors $\CC\lra\CC^{\leq 0}$ (resp. 
$\CC\lra\CC^{\geq 0}$). These functors admit the exact 
left adjoints $\ind_{\fu^{\leq 0}}^{\fu}: \CC^{\leq 0}\lra\CC$ 
(resp. $\ind_{\fu^{\geq 0}}^{\fu}: \CC^{\geq 0}\lra\CC$). 

For example, $\ind_{\fu^{\geq 0}}^{\fu}(M)=\fu\otimes_{\fu^{\geq 0}}M$. 
The triangular decomposition induces an isomorphism of graded  
vector spaces $\ind_{\fu^{\geq 0}}^{\fu}(M)\cong \fu^-\otimes M$. 

\subsection{} For $\lambda\in X$, let us consider on object 
$k^{\lambda}\in\CC^{\geq 0}$ defined as follows. As a graded vector space, 
$k^{\lambda}=k^{\lambda}_{\lambda}=k$. The algebra $\fu^{\geq 0}$ acts on 
$k^{\lambda}$ as follows: $\epsilon_ix=0, K_{\nu}x=\zeta^{\langle\nu,\lambda
\rangle}x\ (i\in I, \nu\in Y, x\in k^{\lambda})$. 

The object $M(\lambda)=\ind_{\fu^{\geq 0}}^{\fu}(k^{\lambda})$ is called 
a {\em (baby) Verma module}. Each $M(\lambda)$ has a unique 
irreducible quotient 
object, to be denoted by $L(\lambda)$. The objects $L(\lambda)\ (\lambda\in X)$ 
are mutually non-isomorphic and every irreducible object in $\CC$ is 
isomorphic to one of them. Note that the category $\CC$ is {\em artinian}, 
i.e. each object of $\CC$ has a finite filtration with irreducible 
quotients.  

For example, $L(0)=\One$.  

\subsection{} Recall that a {\em braiding} on $\CC$ is a collection 
of isomorphisms
$$
R_{M,M'}: M\otimes M'\iso M'\otimes M\ (M,M'\in\CC)
$$
satisfying certain compatibility with the tensor structure (see \cite{kl}IV, 
A.11).  

A {\em balance} on $\CC$ is an automorphism of the identity functor 
$b=\{ b_M: M\iso M\ (M\in\CC)\}$ such that for every $M,N\in\CC$, 
$b_{M\otimes N}\circ(b_M\otimes b_N)^{-1}=R_{N,M}\circ R_{M,N}$ 
(see {\em loc. cit.}). 

\subsection{}
\label{k big} Let $\delta$ denote the determinant of the $I\times I$-matrix 
$(\langle i,j'\rangle)$. From now on we 
assume that $\chara(k)$ does not divide $2\delta$, and   
$k$ contains an element $\zeta'$ such 
that $(\zeta')^{2\delta}=\zeta$; we fix such an element $\zeta'$. 
For a number $q\in\frac{1}{2\delta}\BZ$, $\zeta^q$ will denote 
$(\zeta')^{2\delta q}$. 

\subsection{Theorem} (G. Lusztig) {\em There exists a unique braided 
structure $\{R_{M,N}\}$ on the tensor category $\CC$ such that for any 
$\lambda\in X$ and $M\in\CC$, if $\mu\in X$ is such that   
$M_{\mu'}\neq 0$ implies $\mu'\leq\mu$, then 
$$
R_{L(\lambda),M}(x\otimes y)
=\zeta^{\lambda\cdot\mu}y\otimes x
$$
for any $x\in L(\lambda), y\in M_{\mu}$.} $\Box$ 

\subsection{} Let $n:X\lra\frac{1}{2\delta}\BZ$ be the function defined 
by 
$$
n(\lambda)=\frac{1}{2}\lambda\cdot\lambda-\lambda\cdot\rho_{\ell}.
$$

\subsection{Theorem} {\em There exists a unique balance $b$ on $\CC$ 
such that 
for any $\lambda\in X$, $b_{L(\lambda)}=\zeta^{n(\lambda)}$.} $\Box$ 

\subsection{} The rigid tensor category $\CC$, together with the additional 
structure given by the above braiding and balance, is a {\em ribbon 
category} in the sense of Turaev, cf. \cite{k} and references 
therein.

\section{Braiding local systems}
\label{braiding}  

\subsection{} 
\label{color} For a topological space $T$ and a finite set $J$, $T^J$ 
will denote the space of all maps $J\lra T$ (with the topology 
of the cartesian product). Its points are $J$-tuples 
$(t_j)$ of points of $T$. We denote by $T^{Jo}$ the subspace consisting 
of all $(t_j)$ such that for any $j'\neq j''$ in $J$, $t_{j'}\neq t_{j''}$.   

Let $\nu=\sum \nu_ii\in Y^+$. Let us call an {\em unfolding} of 
$\nu$ a map of finite sets $\pi: J\lra I$ such that $\card(\pi^{-1}(i))=
\nu_i$ for all $i$. Let $\Sigma_{\pi}$ denote the group of all 
automorphisms $\sigma: J\iso J$ such that $\pi\circ\sigma=\pi$. 

The group $\Sigma_{\pi}$ acts on the space $T^J$ in the obvious way. 
We denote by $T^{\nu}$ the quotient space $T^J/\Sigma_{\nu}$. This space 
does not depend, up to a unique isomorphism, on the choice of an unfolding 
$\pi$. The points of $T^{\nu}$ are collections $(t_j)$ of $I$-colored 
points of $T$, such that for any $I$, there are $\nu_i$ points of color $i$. 
We have the canonical projection $T^J\lra T^{\nu}$, also to be denoted 
by $\pi$. We set $T^{\nu o}=\pi(T^{Jo})$. The map $\pi$ restricted to 
$T^{Jo}$ is an unramified Galois covering with the Galois group 
$\Sigma_{\pi}$. 

\subsection{} For a real $r>0$, let $D(r)$ denote the open disk on the 
complex plane $\{ t\in\BC|\ |t|<r\}$ and $\bar{D}(r)$ its closure. 
For $r_1<r_2$, denote by $A(r_1,r_2)$ the open annulus 
$D(r_2)-\bar{D}({r_1})$.  

Set $D=D(1)$. Let $A$ denote the punctured disk $D-\{0\}$. 

\subsection{}   
For an integer $n\geq 1$, consider the space  
$$
E_n=\{ (r_0,\ldots,r_n)\in\BR^{n+1}|\ 0=r_0<r_1<\ldots<r_n=1\}.
$$ 
Obviously, the space $E_n$ is contractible. 

Let $J_1,\ldots,J_n$ be finite sets. Set $J=\coprod_{a=1}^n\ J_a$. 
Note that $D^J=D^{J_1}\times\ldots\times D^{J_n}$.  
Let us untroduce the subspace 
$D^{J_1,\ldots,J_n}\subset E_n\times D^J$. By definition it consists 
of points $((r_a);(t^a_j)\in D^{J_a},\ a=1,\ldots, n)$ such that 
$t^a_j\in A(r_{a-1},r_a)$ for $a=1,\ldots, n$. 
The canonical projection $E_n\times D^J\lra D^J$ induces 
the map
$$ 
m(J_1,\ldots,J_n): D^{J_1,\ldots,J_n}\lra D^J.
$$
The image of the above projection lands in the subspace 
$D^{J_1}\times A^{J_2}\times\ldots\times A^{J_n}\subset D^{J_1}\times
\ldots\times D^{J_n}=D^J$. The induced map 
$$
p(J_1,\ldots,J_n): D^{J_1,\ldots,J_n}\lra D^{J_1}\times A^{J_2}\times\ldots
\times A^{J_n}
$$
is homotopy equivalence. 

Now assume that we have maps $\pi_a: J_a\lra I$ which  
are unfoldings of the elements $\nu_a$. 
Then their sum $\pi: J\lra I$ is an unfolding of $\nu=\nu_1+\ldots+\nu_n$. 
We define the space $D^{\nu_1,\ldots,\nu_n}\subset E_n\times D^{\nu}$ as the 
image of $D^{J_1,\ldots,J_n}$ under the projection 
$\Id\times\pi: E_n\times D^J\lra E_n\times D^{\nu}$. The maps 
$m(J_1,\ldots,J_n)$ and $p(J_1,\ldots,J_n)$ induce 
the maps 
$$
m(\nu_1,\ldots,\nu_n):D^{\nu_1,\ldots,\nu_n}\lra D^{\nu_1+\ldots+\nu_n}
$$
and
$$
p(\nu_1,\ldots,\nu_n):D^{\nu_1,\ldots,\nu_n}\lra D^{\nu_1}\times A^{\nu_2}
\times\ldots\times A^{\nu_n}
$$
respectively, the last map being homotopy equivalence.  

\subsection{} We define the open subspaces  
$$
A^{\nu_1,\ldots,\nu_n}=D^{\nu_1,\ldots,\nu_n}\cap (E_n\times A^{\nu})
$$
and
$$
A^{\nu_1,\ldots,\nu_no}=D^{\nu_1,\ldots,\nu_n}\cap (E_n\times A^{\nu o}).
$$
We have the maps 
$$
m_a(\nu_1,\ldots,\nu_n): D^{\nu_1,\ldots,\nu_n}\lra D^{\nu_1,\ldots,
\nu_{a-1},\nu_a+\nu_{a+1},\nu_{a+2},\ldots,\nu_n}
$$
and
$$
p_a(\nu_1,\ldots,\nu_n): D^{\nu_1,\ldots,\nu_n}\lra 
D^{\nu_1,\ldots,\nu_a}\times A^{\nu_{a+1},\ldots,\nu_n}
$$
$(a=1,\ldots,n-1)$ 
defined in an obvious manner. They induce the maps
$$
m_a(\nu_1,\ldots,\nu_n): A^{\nu_1,\ldots,\nu_n}\lra A^{\nu_1,\ldots,
\nu_{a-1},\nu_a+\nu_{a+1},\nu_{a+2},\ldots,\nu_n}
$$
and
$$
p_a(\nu_1,\ldots,\nu_n): A^{\nu_1,\ldots,\nu_n}\lra 
A^{\nu_1,\ldots,\nu_a}\times A^{\nu_{a+1},\ldots,\nu_n}
$$
and similar maps between "$o$"-ed spaces. All the maps $p$ 
are homotopy equivalences. 

All these maps satisfy some obvious compatibilities. We will need the following 
particular case. 

\subsection{} 
\label{rhomb} The {\em rhomb} diagram below commutes. 

\begin{center}
  \begin{picture}(14,8)
    \put(7,8){\makebox(0,0){$D^{\nu_1+\nu_2+\nu_3}$}}
    \put(5,6){\makebox(0,0){$D^{\nu_1,\nu_2+\nu_3}$}}
    \put(9,6){\makebox(0,0){$D^{\nu_1+\nu_2,\nu_3}$}}
    \put(3,4){\makebox(0,0){$D^{\nu_1}\times A^{\nu_1+\nu_2}$}}
    \put(7,4){\makebox(0,0){$D^{\nu_1,\nu_2,\nu_3}$}}
    \put(11,4){\makebox(0,0){$D^{\nu_1+\nu_2}\times A^{\nu_3}$}}
    \put(5,2){\makebox(0,0){$D^{\nu_1}\times A^{\nu_1,\nu_2}$}}
    \put(9,2){\makebox(0,0){$D^{\nu_1,\nu_2}\times A^{\nu_3}$}}
    \put(7,0){\makebox(0,0){$D^{\nu_1}\times A^{\nu_2}\times A^{\nu_3}$}}

    \put(5.5,6.5){\vector(1,1){1}}
    \put(8.5,6.5){\vector(-1,1){1}}
    \put(4.5,5.5){\vector(-1,-1){1}}
    \put(6.5,4.5){\vector(-1,1){1}}
    \put(7.5,4.5){\vector(1,1){1}}
    \put(9.5,5.5){\vector(1,-1){1}} 
    \put(4.5,2.5){\vector(-1,1){1}}
    \put(6.5,3.5){\vector(-1,-1){1}}
    \put(7.5,3.5){\vector(1,-1){1}}
    \put(9.5,2.5){\vector(1,1){1}}
    \put(5.5,1.5){\vector(1,-1){1}}
    \put(8.5,1.5){\vector(-1,-1){1}}

    \put(5.5,7){\makebox(0,0){$m$}}
    \put(8.5,7){\makebox(0,0){$m$}}
    \put(3.5,5){\makebox(0,0){$p$}}
    \put(5.5,5){\makebox(0,0){$m$}}
    \put(8.5,5){\makebox(0,0){$m$}}
    \put(10.5,5){\makebox(0,0){$p$}}
    \put(3.5,3){\makebox(0,0){$m$}}
    \put(5.5,3){\makebox(0,0){$p$}}
    \put(8.5,3){\makebox(0,0){$p$}}
    \put(10.5,3){\makebox(0,0){$m$}}
    \put(5.5,1){\makebox(0,0){$p$}}
    \put(8.5,1){\makebox(0,0){$p$}}
    
  \end{picture} 
\end{center}

\subsection{} We will denote by $\CA^o$ and call an {\em open 
$I$-coloured configuration space} the collection of all spaces 
$\{ A^{\nu_1,\ldots,\nu_n o}\}$ together with the maps 
$\{ m_a(\nu_1,\ldots,\nu_n), 
p_a(\nu_1,\ldots,\nu_n)\}$ between their various products. 

We will call a {\em local system} over $\CA^o$, or a {\em braiding local 
system} a collection of data (a), (b) below satisfying the property 
(c) below. 

(a) A local system $\CI_{\mu}^{\nu}$ over $A^{\nu o}$ given for any 
$\nu\in  Y^+, \mu\in X$. 

(b) An isomorphism 
$\phi_{\mu}(\nu_1,\nu_2): m^*\CI_{\mu}^{\nu_1+\nu_2}\iso p^*(\CI_{\mu}^{\nu_1}
\boxtimes\CI_{\mu-\nu_1}^{\nu_2})$
given for any $\nu_1,\nu_2\in  Y^+, \mu\in X$. 

Here $p=p(\nu_1,\nu_2)$ and 
$m=m(\nu_1,\nu_2)$ are the arrows in the diagram 
$A^{\nu_1 o}\times A^{\nu_2 o}\overset{p}{\lla} A^{\nu_1,\nu_2 o}
\overset{m}{\lra} A^{\nu_1+\nu_2 o}$. 

The isomorphisms $\phi_{\mu}(\nu_1,\nu_2)$ are called the {\em factorization 
isomorphisms}. 

(c) (The {\em associativity} of factorization isomorphisms.) For any  
$\nu_1, \nu_2, \nu_3\in  Y^+, \mu\in X$, the octagon below commutes. 
Here the maps 
$m, p$ are the maps in the rhombic diagram above, with $D$, $A$ replaced 
by $A^o$.   

\begin{center}
  \begin{picture}(14,12)
    \put(4,12){\makebox(0,0){$m^*m_1^*\CI^{\nu_1+\nu_2+\nu_3}_{\mu}$}}
    \put(10,12){\makebox(0,0){$m^*m_2^*\CI^{\nu_1+\nu_2+\nu_3}_{\mu}$}}
    \put(0,8){\makebox(0,0){$m^*p^*(\CI^{\nu_1}_{\mu}\boxtimes\CI^{\nu_2+\nu_3}
    _{\mu-\nu_1})$}}
    \put(14,8){\makebox(0,0){$m^*p^*(\CI^{\nu_1+\nu_2}_{\mu}\boxtimes\CI^{\nu_3}
    _{\mu-\nu_1-\nu_2})$}}
    \put(0,4){\makebox(0,0){$p^*m^*(\CI^{\nu_1}_{\mu}\boxtimes\CI^{\nu_2+\nu_3}
    _{\mu-\nu_1})$}}
    \put(14,4){\makebox(0,0){$p^*m^*(\CI^{\nu_1+\nu_2}_{\mu}\boxtimes\CI^{\nu_3}
    _{\mu-\nu_1-\nu_2})$}}
    \put(3,0){\makebox(0,0){$p_1^*p^*(\CI_{\mu}^{\nu_1}\boxtimes\CI_{\mu-\nu_1}
    ^{\nu_2}\boxtimes\CI_{\mu-\nu_1-\nu_2}^{\nu_3})$}}
    \put(11,0){\makebox(0,0){$p_2^*p^*(\CI_{\mu}^{\nu_1}\boxtimes\CI_{\mu-\nu_1}
    ^{\nu_2}\boxtimes\CI_{\mu-\nu_1-\nu_2}^{\nu_3})$}}

    \put(6,12){\line(1,0){2}}
    \put(6,11.9){\line(1,0){2}}
    \put(3.5,11.5){\vector(-1,-1){3}}
    \put(10,11.5){\vector(1,-1){3}}
    \put(0,7.5){\line(0,-1){3}}
    \put(.1,7.5){\line(0,-1){3}}
    \put(14,7.5){\line(0,-1){3}}
    \put(13.9,7.5){\line(0,-1){3}}
    \put(.5,3.5){\vector(1,-1){3}}
    \put(13.5,3.5){\vector(-1,-1){3}}
    \put(6.4,0){\line(1,0){1}}
    \put(6.4,.1){\line(1,0){1}}

    \put(0,10){\makebox(0,0){$\phi_{\mu}(\nu_1,\nu_2+\nu_3)$}}
    \put(14,10){\makebox(0,0){$\phi_{\mu}(\nu_1+\nu_2,\nu_3)$}}
    \put(0,2){\makebox(0,0){$\phi_{\mu-\nu_1}(\nu_2,\nu_3)$}}
    \put(14,2){\makebox(0,0){$\phi_{\mu}(\nu_1,\nu_2)$}}

  \end{picture} 
\end{center}

Written more concisely, the axiom (c) reads as a "cocycle" condition 

(c)$'$  $\phi_{\mu}(\nu_1,\nu_2)\circ\phi_{\mu}(\nu_1+\nu_2,\nu_3)=
\phi_{\mu-\nu_1}(\nu_2,\nu_3)\circ\phi_{\mu}(\nu_1,\nu_2+\nu_3)$. 

\subsection{Correctional lemma} {\em Assume we are given the data 
{\em (a), (b)} as above, with  
one-dimensional local systems $\CI_{\mu}^{\nu}$. 
Then there exist a collection of constants $c_{\mu}(\nu_1,\nu_2)\in k^*\ 
(\mu\in X,\nu_1,\nu_2\in  Y^+)$ such that the corrected isomorphisms 
$\phi'_{\mu}(\nu_1,\nu_2)=c_{\mu}(\nu_1,\nu_2)\phi_{\mu}(\nu_1,\nu_2)$ 
satisfy the associativity axiom {\em (c)}}. $\Box$  

\subsection{} The notion of a morphism between two braiding local 
systems is defined in an obvious way. This defines the category 
of braiding local systems, $\Bls$.  

Suppose that $\CI=\{\CI_{\mu}^{\nu}; \phi_{\mu}(\nu_1,\nu_2)\}$ and 
$\CJ=\{\CJ_{\mu}^{\nu}; \psi_{\mu}(\nu_1,\nu_2)\}$ are two braiding local 
systems. Their 
{\em tensor product} $\CI\otimes\CJ$ is defined by 
$(\CI\otimes\CJ)_{\mu}^{\nu}=\CI_{\mu}^{\nu}\otimes\CJ_{\mu}^{\nu}$, 
the factorization isomorphisms being $\phi_{\mu}(\nu_1,\nu_2)\otimes 
\psi_{\mu}(\nu_1,\nu_2)$. This makes $\Bls$ a tensor category. 
The subcategory of one-dimensional braiding local systems is a Picard  
category. 

\subsection{Example} 
\label{sign} {\em Sign local system.} Let $L$ be a one-dimensional 
vector space. Let $\nu\in Y^+$. Let $\pi: J\lra I$ be an unfolding of $\nu$, 
whence the canonical projection $\pi: A^{J\circ}\lra A^{\nu\circ}$.  
Pick a point $\tilde{x}=(x_j)\in A^{J\circ}$ with all $x_j$ being {\em 
real}; let $x=\pi(\tilde{x})$. The choice of base points defines the 
homomorphism $\pi_1(A^{\nu\circ};x)\lra\Sigma_{\pi}$. Consider its 
composition with the sign map $\Sigma_{\pi}\lra\mu_2\hra k^*$. Here 
$\mu_2=\{\pm 1\}$ is the group of square roots of $1$ in $k^*$. 
We get a map $s: \pi_1(A^{\nu\circ};x)\lra k^*$. It defines a local system 
$\Sign^{\nu}(L)$ over $A^{\nu\circ}$ whose stalk at $x$ is equal to $L$.  
This local system does not depend (up to the 
unique isomorphism) on the choices made. 

Given a diagram
$$
A^{\nu_1\circ}\times A^{\nu_2\circ}\overset{p}{\lla} A^{\nu_1,\nu_2\circ}
\overset{m}{\lra}A^{\nu_2+\nu_2\circ}
$$
we can choose base point $x_i\in A^{\nu_i}$ and $x\in A^{\nu_1+\nu_2\circ}$ 
such that there exists $y\in A^{\nu_1,\nu_2}$ with $p(y)=(x_1,x_2)$ and 
$m(y)=x$. Let 
$$
\phi^{L_1,L_2}(\nu_1,\nu_2): m^*\Sign^{\nu_1+\nu_2}(L_1\otimes L_2)\iso 
p^*(\Sign^{\nu_1}(L_1)\boxtimes\Sign^{\nu_2}(L_2))
$$
be the unique isomorphism equal to $\Id_{L_1\otimes L_2}$ over $y$. 
These isomorphisms satisfy the associativity condition (cf. (c)$'$ above). 

Let us define the brading local system $\Sign$ by $\Sign_{\mu}^{\nu}=
\Sign^{\nu}(\One)$, $\phi_{\mu}(\nu_1,\nu_2)=\phi^{\One,\One}(\nu_1,\nu_2)$. 
Here $\One$ denotes the standard vector space $k$. 

\subsection{Example} 
\label{stand loc} {\em Standard local systems $\CJ$, $\CI$.} Let 
$\nu\in Y^+$, let $\pi: J\lra I$ be an unfolding of $\nu$, $n=\card(J)$, 
$\pi: A^{J\circ}\lra A^{\nu\circ}$ the 
corresponding projection. Let $L$ be a one-dimensional vector space. 
For each isomorphism $\sigma: J\iso \{1,\ldots, n\}$, choose a point $x(\sigma)= 
(x_j)\in A^{J\circ}$ such that all $x_j$ are real and positive and for each 
$j', j''$ such that $\sigma(j')<\sigma(j'')$, we have $x_{j'}<x_{j''}$.
  
Let us define a local system $\CJ_{\mu}^{\pi}(L)$ over $A^{J\circ}$ as follows.  
We set $\CJ_{\mu}^{\pi}(L)_{x(\sigma)}=L$ for all $\sigma$. We define the 
monodromies 
$$
T_{\gamma}: \CJ_{\mu}^{\pi}(L)_{x(\sigma)}\lra \CJ_{\mu}^{\pi}(L)_{x(\sigma')}
$$
along the homotopy classes of paths $\gamma: x_{\sigma}\lra x_{\sigma'}$ 
which generate the fundamental groupoid. Namely, let $x_{j'}$ and $x_{j''}$ 
be some neighbour points in $x(\sigma)$, with $x_{j'}<x_{j''}$. 
Let $\gamma(j',j'')^+$ (resp. $\gamma(j',j'')^-$) be the paths corresponding 
to the movement of $x_{j''}$ in the upper (resp. lower) 
hyperplane to the left from $x_{j'}$ position. We set 
$$
T_{\gamma(j',j'')^{\pm}}=\zeta^{\mp\pi(j')\cdot\pi(j'')}.
$$
Let $x_j$ be a point in $x(\sigma)$ closest to $0$. Let $\gamma(j)$ be the path 
corresponding to the counterclockwise travel of $x_j$ around $0$. We set 
$$
T_{\gamma(j)}=\zeta^{2\mu\cdot\pi(j)}.
$$
The point is that the above fromulas give a well defined morphism from 
the fundamental groupoid to the groupoid of one-dimensional vector spaces. 

The local system $\CJ_{\mu}^{\pi}(L)$ admits an obvious 
$\Sigma_{\pi}$-equivariant structure. This defines the local system 
$\CJ_{\mu}^{\nu}(L)$ over $A^{\nu\circ}$. Given a diagram 
$$
A^{\nu_1\circ}\times A^{\nu_2\circ}\overset{p}{\lla} A^{\nu_1+\nu_2\circ}
\overset{m}{\lra} A^{\nu_1+\nu_2\circ},
$$
let 
$$
\phi_{\mu}^{L_1,L_2}(\nu_1,\nu_2): m^*\CJ^{\nu_1+\nu_2}_{\mu}(L_1\otimes L_2)
\iso p^*(\CJ^{\nu_1}_{\mu}(L_1)\boxtimes\CJ^{\nu_2}_{\mu-\nu_1}(L_2))
$$
be the unique isomorphism equal to $Id_{L_1\otimes L_2}$ over compatible 
base points. Here the compatibility is understood in the same sense as in the 
previous subsection. These isomorphisms satisfy the associativity 
condition.  

We define the braiding local system $\CJ$ by $\CJ_{\mu}^{\nu}=\CJ_{\mu}^{\nu}
(\One)$, $\phi_{\mu}(\nu_1,\nu_2)=\phi_{\mu}^{\One,\One}(\nu_1,\nu_2)$. 

We define the braiding local system $\CI$ by $\CI=\CJ\otimes\Sign$. The local 
systems $\CJ$, $\CI$ are called the standard local systems. In the sequel 
we will mostly need the local system $\CI$.

\section{Factorizable sheaves}
\label{factoriz}

\subsection{} In the sequel for each $\nu\in Y^+$, we will denote by 
$\CS$ the stratification on the space $D^{\nu}$ the closures of whose 
strata are various intersections of hypersurfaces given by the equations 
$t_j=0,\ t_{j'}=t_{j''}$. The same letter will denote the induced 
stratifications on its subspaces.  

\subsection{} Set  
$\CI_{\mu}^{\nu\bullet}=j_{!*}\CI_{\mu}^{\nu}[\dim A^{\nu\circ}]
\in\CM(A^{\nu};\CS)$ 
where $j: A^{\nu\circ}\hra A^{\nu}$ is the embedding. 
   
Given a diagram 
$$
A^{\nu_1}\times A^{\nu_2}\overset{p}{\lla} A^{\nu_1,\nu_2}\overset{m}
{\lra} A^{\nu_1+\nu_2},
$$
the structure isomorphisms $\phi_{\mu}(\nu_1,\nu_2)$ of the local 
system $\CI$ induce the isomorphisms, to be denoted by the same 
letters, 
$$
\phi_{\mu}(\nu_1,\nu_2): m^*\CI^{\nu_1+\nu_2\bullet}\iso 
p^*(\CI^{\nu_1\bullet}\boxtimes\CI^{\nu_2\bullet}).
$$
Obviously, these isomorphisms satisfy 
the associativity axiom. 

\subsection{} Let us fix a coset $c\in X/Y$. We will regard $c$ as a 
subset of $X$.   
We will call a {\em factorizable sheaf} $\CM$ {\em supported at $c$} 
a collection 
of data (w), (a), (b) below satisfying the axiom (c) below. 

(w) An element $\lambda=\lambda(\CM)\in c$. 

(a) A perverse sheaf $\CM^{\nu}\in\CM(D^{\nu};\CS)$ given for each $\nu\in Y^+$. 

(b) An isomorphism $\psi(\nu_1,\nu_2): m^*\CM^{\nu_1+\nu_2}\iso 
p^*(\CM^{\nu_1}\boxtimes\CI_{\lambda-\nu_1}^{\nu_2\bullet})$ 
given for any $\nu_1, \nu_2\in Y^+$. 

Here $p, m$ denote the arrows in the diagram $D^{\nu_1}\times A^{\nu_2}
\overset{p}{\lla} D^{\nu_1,\nu_2}\overset{m}{\lra} D^{\nu_1+\nu_2}$. 

The isomorphisms $\psi(\nu_1,\nu_2)$ are called the {\em factorization 
isomorphisms}. 

(c) For any $\nu_1, \nu_2, \nu_3\in Y^+$, the following  
{\em associativity condition} is fulfilled:  
$$
\psi(\nu_1,\nu_2)\circ\psi(\nu_1+\nu_2,\nu_3)=
\phi_{\lambda-\nu_1}(\nu_2,\nu_3)\circ\psi(\nu_1,\nu_2+\nu_3).
$$

We leave to the reader to draw the whole octagon expressing this axiom. 

\subsection{} Let $\CM=\{\CM^{\nu};\ \psi(\nu_1,\nu_2)\}$ be a factorizable 
sheaf supported at a coset $c\in X/Y$, $\lambda=\lambda(\CM)$. For each 
$\lambda'\geq\lambda, \nu\in Y^+$, define a sheaf $\CM_{\lambda'}^{\nu}\in
\CM(D^{\nu};\CS)$ by  
$$
\CM_{\lambda'}^{\nu}=\left\{ \begin{array}{ll}\iota(\lambda'-\lambda)_*
\CM^{\nu-\lambda'+\lambda}&\mbox{ if }\nu-\lambda'+\lambda\in Y^+\\
0&\mbox{ otherwise.}\end{array}\right.
$$
Here 
$$
\iota(\nu'): D^{\nu}\lra D^{\nu+\nu'}
$$
denotes the closed embedding adding $\nu'$ points sitting at the origin. 
The factorization isomorphisms $\psi(\nu_1,\nu_2)$ induce  
similar isomorphisms 
$$
\psi_{\lambda'}(\nu_1,\nu_2): m^*\CM_{\lambda'}(\nu_1+\nu_2)\iso 
p^*(\CM_{\lambda'}^{\nu_1}\boxtimes\CI_{\lambda'-\nu_1}^{\nu_2})\ (\lambda'
\geq\lambda)
$$

\subsection{} Let $\CM, \CN$ be two factorizable sheaves supported 
at $c$. Let $\lambda\in X$ 
be such that $\lambda\geq\lambda(\CM)$ and $\lambda\geq\lambda(\CN)$. 
For $\nu\geq\nu'$ in $Y^+$, consider the following composition
$$
\tau_{\lambda}(\nu,\nu'):\Hom(\CM_{\lambda}^{\nu},\CN_{\lambda}^{\nu})
\overset{m_*}{\lra}\Hom(m^*\CM_{\lambda}^{\nu},m^*\CN_{\lambda}^{\nu})
\overset{\psi(\nu',\nu-\nu')_*}{\iso}
$$
$$
\iso\Hom(p^*(\CM_{\lambda}^{\nu'}
\boxtimes\CI_{\lambda-\nu'}^{\nu-\nu'\bullet}),p^*(\CN_{\lambda}^{\nu'}
\boxtimes\CI_{\lambda-\nu'}^{\nu-\nu'\bullet}))=\Hom(\CM_{\lambda}^{\nu'},
\CN_{\lambda}^{\nu'}).
$$
Let us define the space of homomorphisms $\Hom(\CM,\CN)$ by 
$$
\Hom(\CM,\CN)=\dirlim_{\lambda}\invlim_{\nu}\Hom(\CM_{\lambda}^{\nu},
\CN_{\lambda}^{\nu})
$$
Here the inverse limit is taken over $\nu\in Y^+$, the transition maps being 
$\tau_{\lambda}(\nu,\nu')$ and the direct limit is taken over $\lambda\in X$ 
such that $\lambda\geq\lambda(\CM), \lambda\geq\lambda(\CN)$,  
the transition maps being induced by the obvious isomorphisms
$$
\Hom(\CM_{\lambda}^{\nu},\CN_{\lambda}^{\nu})=
\Hom(\CM_{\lambda+\nu'}^{\nu+\nu'},\CN_{\lambda+\nu'}^{\nu+\nu'})\ 
(\nu'\in Y^+). 
$$
With these spaces of homomorphisms and the obvious compositions, 
the factorizable sheaves supported at $c$ form the category, to be 
denoted by $\tFS_c$. By definition, the category $\tFS$ is the direct 
product $\prod_{c\in X/Y}\ \tFS_c$. 

\subsection{Finite sheaves} Let us call a factorizable sheaf 
$\CM=\{\CM^{\nu}\}\in\tFS_c$ {\em finite} if there exists only a finite 
number of $\nu\in Y^+$ such that the conormal bundle of the origin 
$O\in\CA^{\nu}$ is contained in the singular support of $\CM^{\nu}$. 
Let $\FS_c\subset\tFS_c$ be the full subcategory of finite 
factorizable sheaves. We define the category $\FS$ by 
$\FS=\prod_{c\in X/Y}\FS_c$. One proves (using the lemma below) that 
$\FS$ is an abelian category.  

\subsection{Stabilization Lemma} {\em Let $\CM, \CN\in\FS_c,\ 
\mu\in X_c, \mu\geq\lambda(\CM), \mu\geq\lambda(\CN)$. There exists 
$\nu_0\in Y^+$ such that for all $\nu\geq\nu_0$ the canonical maps 
$$
\Hom(\CM,\CN)\lra\Hom(\CM_{\mu}^{\nu},\CN_{\mu}^{\nu})
$$
are isomorphisms.} $\Box$   

\subsection{Standard sheaves} 
\label{stand} Given $\mu\in X$, let us define the 
"standard sheaves" $\CM(\mu), \CD\CM(\mu)$ and $\CL(\mu)$ 
supported at the coset $\mu+Y$, by 
$\lambda(\CM(\mu))=\lambda(\CD\CM(\mu))=\lambda(\CL(\mu))=\mu$;  
$$
\CM(\mu)^{\nu}=j_!\CI^{\nu\bullet}_{\mu};\ 
\CD\CM(\mu)^{\nu}=j_*\CI^{\nu\bullet}_{\mu};\ 
\CL(\mu)^{\nu}=j_{!*}\CI^{\nu\bullet}_{\mu},
$$
$j$ being the embedding $A^{\nu}\hra D^{\nu}$. The factorization maps 
are defined by functoriality from the similar maps for $\CI^{\bullet}$. 

One proves that all these sheaves are finite.

\section{Tensor product}

In this section we will give (a sketch of) the construction of the 
tensor structure on the category $\tFS$. We will make the assumption of 
~\ref{k big}\footnote{Note that these assumptions are not necessary for the 
construction of the tensor structure. They are essential, however, 
for the construction of braiding.}. 

\subsection{} For $z\in\BC$ and a real positive $r$, let $D(z;r)$ 
denote the open disk $\{ t\in\BC|\ |t-z|<r\}$ and $\bar{D}(z;r)$ its 
closure. 

\subsection{} 
\label{spaces} For $\nu\in Y^+$, let us define the space $D^{\nu}(2)$ 
as the product $A\times D^{\nu}$. Its points will be denoted 
$(z;(t_j))$ where $z\in A, (t_j)\in D^{\nu}$. Let us define the open 
subspaces 
$$
A^{\nu}(2)=\{(z;(t_j))\in D^{\nu}(2)|\ t_j\neq 0, z\mbox{ for all }j\};\ 
A^{\nu}(2)^{\circ}=A^{\nu}(2)\cap(A\times A^{\nu\circ}). 
$$
For $\nu,\nu'\in Y^+$, let us define the space 
$D^{\nu,\nu'}(2)$ as the subspace of $\BR_{>0}\times D^{\nu+\nu'}(2)$ 
consisting of all elements $(r;z;(t_j))$ such that $|z|<r<1$ and $\nu$ 
of the points $t_j$ live inside the disk $D(r)$ and $\nu'$ of them inside  
the annulus $A(r,1)$. 

We have a diagram

(a) $D^{\nu}(2)\times A^{\nu'}\overset{p}{\lla} D^{\nu,\nu'}(2)
\overset{m}{\lra} D^{\nu+\nu'}(2)$. 

Here $p((r;z;(t_j)))=((z;(t_{j'})), (t_{j''}))$ where $t_{j'}$ (resp. $t_{j''}$) 
being the points from the collection $(t_j)$ lying in $D(r)$ 
(resp. in $A(r,1)$); 
$m((r;z;(t_j)))=(z;(t_j))$. The map $p$ is a homotopy equivalence. 

For $\nu_1,\nu_2,\nu\in Y^+$, let $D^{\nu_1;\nu_2;\nu}(2)$ be the 
subspace of $\BR_{>0}\times\BR_{>0}\times D^{\nu_1+\nu_2+\nu}(2)$ 
consisting of all elements $(r_1;r_2;z;(t_j))$ such that 
$\bar{D}(r_1)\cup\bar{D}(z;r_2)\subset D$; $\bar{D}(r_1)\cap\bar{D}(z;r_2)
=\emptyset$; $\nu_1$ of the points $(t_j)$ lie inside $D(r_1)$, 
$\nu_2$ of them lie inside $D(z;r_2)$ and $\nu$ of them lie 
inside $D-(\bar{D}(r_1)\cup\bar{D}(z;r_2))$. 

We have a diagram 

(b) $D^{\nu_1}\times D^{\nu_2}\times A^{\nu}(2)\overset{p}{\lla} 
D^{\nu_1;\nu_2;\nu}(2)\overset{m}{\lra} D^{\nu_1+\nu_2+\nu}(2)$. 

Here $p((r_1;r_2;z;(t_j)))=((t_{j'}); (t_{j''}-z); (t_{j'''}))$ where 
$t_{j'}$ (resp. $t_{j''}, t_{j'''}$) are the points lying inside $D(r_1)$ 
(resp. $D(z;r_2), D-(\bar{D}(r_1)\cup\bar{D}(z;r_2)))$; 
$m((r_1;r_2;z;(t_j)))=(z;(t_j))$. The map $p$ is a homotopy equivalence. 

\subsection{}    
We set $A^{\nu,\nu'}(2)=D^{\nu,\nu'}(2)\cap(\BR_{>0}
\times A^{\nu+\nu'}(2))$; 
$A^{\nu,\nu^{\prime}}(2)^{\circ}=D^{\nu,\nu'}(2)\cap(\BR_{>0}\times 
A^{\nu+\nu'}(2)^{\circ})$; $A^{\nu_1;\nu_2;\nu}(2)=D^{\nu_1;\nu_2;\nu}(2)
\cap(\BR_{>0}\times\BR_{>0}\times A^{\mu_1+\nu_2+\nu}(2))$; 
$A^{\nu_1;\nu_2;\nu}(2)^{\circ}=D^{\nu_1;\nu_2;\nu}(2)\cap(\BR_{>0}
\times\BR_{>0}\times A^{\nu_1+\nu_2+\nu}(2)^{\circ})$. 

\subsection{} Given $\mu_1,\mu_2\in X, \nu\in Y^+$, choose an unfolding 
$\pi: J\lra I$ of the element $\nu$.  
In the same manner as in \ref{stand loc}, we define the   
one-dimensional local system $\CJ_{\mu_1,\mu_2}^{\nu}$ over 
$A^{\nu}(2)^{\circ}$ with the following monodromies: the monodromy 
around a loop corresponding to the counterclockwise travel of the point 
$z$ around $0$ (resp. $t_j$ around $0$, $t_j$ around $z$, 
$t_{j'}$ around $t_{j''}$) is equal to the multiplication 
by $\zeta^{-2\mu_1\cdot\mu_2}$ (resp. $\zeta^{2\mu_1\cdot\pi(j)}$, 
$\zeta^{2\mu_2\cdot\pi(j)}$, $\zeta^{-2\pi(j')\cdot\pi(j'')}$).  

As in {\em loc. cit.}, one defines isomorphisms

(a) $\phi_{\mu_1,\mu_2}(\nu,\nu'): m^*\CJ_{\mu_1,\mu_2}^{\nu+\nu'}\iso 
p^*(\CJ_{\mu_1,\mu_2}^{\nu}\boxtimes\CJ_{\mu_1+\mu_2-\nu}^{\nu'})$ 

where $p$, $m$ are the morphisms in the diagram \ref{spaces} (a)  
(restricted to the $A^{\circ}$-spaces) and  

(b) $\phi_{\mu_1,\mu_2}(\nu_1;\nu_2;\nu): m^*\CJ_{\mu_1,\mu_2}^{\nu_1+
\nu_2+\nu}\iso p^*(\CJ_{\mu_1}^{\nu_1}\boxtimes\CJ_{\mu_2}^{\nu_2}\boxtimes
\CJ_{\mu_1-\nu_1,\mu_2-\nu_2}^{\nu})$ 

where $p, m$ are the morphisms in the diagram \ref{spaces} (b) (restricted  
to the $A^{\circ}$-spaces),   
which satisfy the cocycle conditions 

(c) $\phi_{\mu_1,\mu_2}(\nu,\nu')\circ\phi_{\mu_1,\mu_2}(\nu+\nu',\nu'')=
\phi_{\mu_1+\mu_2-\nu}(\nu',\nu'')\circ\phi_{\mu_1,\mu_2}(\nu,\nu'+\nu'')$ 

(d) $(\phi_{\mu_1}(\nu_1,\nu'_1)\boxtimes\phi_{\mu_2}(\nu_2,\nu'_2))
\circ\phi_{\mu_1,\mu_2}(\nu_1+\nu'_1;\nu_2+\nu'_2;\nu)=$

$=\phi_{\mu_1-\nu_1,\mu_2-\nu_2}(\nu'_1;\nu'_2;\nu)\circ\phi_{\mu_1,\mu_2}
(\nu_1;\nu_2;\nu+\nu'_1+\nu'_2)$  

(we leave to the reader the definition of the corresponding spaces).  

\subsection{} Let us consider the sign local systems introduced in \ref{sign}. 
We will keep the same notation $\Sign^{\nu}$ for the inverse image 
of the local system $\Sign^{\nu}$ under the forgetting of 
$z$ map $A^{\nu}(2)^{\circ}\lra A^{\nu\circ}$. We have the factorization 
isomorphisms

(a) $\phi^{\Sign}(\nu,\nu'): m^*\Sign^{\nu+\nu'}\iso p^*(\Sign^{\nu}\boxtimes
\Sign^{\nu'})$; 

(b) $\phi^{\Sign}(\nu_1;\nu_2;\nu): m^*\Sign^{\nu_1+\nu_2+\nu}\iso p^*(
\Sign^{\nu_1}\boxtimes\Sign^{\nu_2}\boxtimes\Sign^{\nu})$ 

which satisfy the cocycle conditions similar to (c), (d) above.    

\subsection{}  
We define the local systems $\CI^{\nu}_{\mu_1,\mu_2}$ over the spaces 
$A^{\nu}(2)^{\circ}$ by  

$\CI_{\mu_1,\mu_2}^{\nu}=\CJ_{\mu_1,\mu_2}^{\nu}\otimes\Sign^{\nu}$.  
  
The collection of local systems $\{\CI_{\mu_1,\mu_2}^{\nu}\}$ together 
with the maps  
$\phi^{\CI}_{\mu_1,\mu_2}(\nu,\nu')=\phi^{\CJ}_{\mu_1,\mu_2}(\nu,\nu')\otimes
\phi^{\Sign}(\nu,\nu')$ and    
$\phi^{\CI}_{\mu_1,\mu_2}(\nu_1;\nu_2;\nu)=\phi^{\CJ}_{\mu_1,\mu_2}(\nu_1;\nu_2;
\nu)\otimes\phi^{\Sign}(\nu_1;\nu_2;\nu)$, forms an object $\CI(2)$ which 
we call 
a {\em standard braiding local system over the configuration space $\CA(2)^{\circ}
=\{ A^{\nu}(2)^{\circ}\}$}. It is unique up to a (non unique) isomorphism.  
We fix such a local system.  

\subsection{} We set $\CI_{\mu_1,\mu_2}^{\nu\bullet}=
j_{!*}\CI_{\mu_1,\mu_2}^{\nu}
[\dim A^{\nu}(2)^{\circ}]$ where $j: A^{\nu}(2)^{\circ}\hra A^{\nu}(2)$ 
is the open embedding. It is an object of the category $\CM(A^{\nu}(2);\CS)$ 
where $\CS$ is the evident stratification. The factorization isomorphisms 
for the local system $\CI$ induce the analogous isomorphisms 
between these sheaves, to be denoted by the same letter. The collection 
of these sheaves and factorization isomorphisms will be denoted 
$\CI(2)^{\bullet}$.     

\subsection{} 
Suppose we are given two factorizable sheaves $\CM, \CN$. Let us call their 
{\em gluing}, and denote by $\CM\boxtimes\CN$, the collection of 
perverse sheaves $(\CM\boxtimes\CN)^{\nu}$ over the spaces 
$D^{\nu}(2)$ $(\nu\in Y^+)$ together with isomorphisms
$$
\psi(\nu_1;\nu_2;\nu): m^*(\CM\boxtimes\CN)^{\nu}\iso p^*(\CM^{\nu_1}
\boxtimes\CN^{\nu_2}\boxtimes\CI^{\nu\bullet}_{\lambda(\CM)-\nu_1,
\lambda(\CN)-\nu_2}),
$$
$p, m$ being the maps in the diagram \ref{spaces} (b), which satisfy the 
cocycle condition

$(\psi^{\CM}(\nu_1,\nu'_1)\boxtimes\psi^{\CN}(\nu_2,\nu'_2))\circ
\psi(\nu_1+\nu'_1;\nu_2+\nu'_2;\nu)=$

$=\phi_{\lambda(\CM)-\nu_1,\lambda(\CN)-\nu_2}(\nu'_1;\nu'_2;\nu)\circ
\psi(\nu_1;\nu_2;\nu+\nu_1'+\nu_2')$ 

for all $\nu_1,\nu_1',\nu_2,\nu_2',\nu\in Y^+$. 

Such a gluing exists and is unique, up to a unique isomorphism. The 
factorization isomorphisms $\phi_{\mu_1,\mu_2}(\nu_1;\nu_2;\nu)$ 
for $\CI(2)^{\bullet}$ and the ones for $\CM, \CN$,  
induce the isomorphisms 
$$
\psi^{\CM\boxtimes\CN}(\nu,\nu'): m^*(\CM\boxtimes\CN)^{\nu+\nu'}\iso 
p^*((\CM\boxtimes\CN)^{\nu}\boxtimes\CI^{\nu'\bullet}_{\lambda(\CM)+
\lambda(\CN)-\nu}).
$$
satisfying the obvious cocycle condition. 

\subsection{} 
\label{def tens} Now we can define the tensor product $\CM\otimes\CN\in\tFS$. 
Namely, set $\lambda(\CM\otimes\CN)=\lambda(\CM)+\lambda(\CN)$. For 
each $\nu\in Y^+$, set 
$$
(\CM\otimes\CN)^{\nu}=\Psi_{z\ra 0}((\CM\boxtimes\CN)^{\nu}).
$$
Here $\Psi_{z\ra 0}: \CM(D^{\nu}(2))\lra\CM(D^{\nu})$ denotes the functor 
of nearby cycles for the function $D^{\nu}(2)\lra D$ sending $(z;(t_j))$ to 
$z$. Note that 
$$
\Psi_{z\ra 0}(\CI_{\mu_1,\mu_2}^{\nu})=\CI_{\mu_1+\mu_2}^{\nu}.
$$
The factorization isomorphisms $\psi^{\CM\boxtimes\CN}$ induce the factorization 
isomorphisms between the sheaves $(\CM\otimes\CN)^{\nu}$. 
This defines a factorizable sheaf $\CM\otimes\CN$. 

One sees at once that this construction is functorial; thus it 
defines a functor of tensor product $\otimes: \tFS\times\tFS\lra\tFS$. 

The subcategory $\FS\subset\tFS$ is stable under the tensor product. 
The functor $\otimes:\ \FS\times\FS\lra\FS$ extends uniquely to a 
functor $\otimes:\ \FS\otimes\FS\lra\FS$ (for the discussion of the 
tensor product of abelian categories, see \cite{d2} 5).    

\subsection{} The half-circle travel of the point $z$ around $0$ from $1$ to 
$-1$ in the upper halfplane defines the braiding isomorphisms 
$$
{R}_{\CM,\CN}: \CM\otimes\CN\iso\CN\otimes\CM.
$$
We will not describe here the precise definition of the associativity 
isomorphisms for the tensor product $\otimes$. We just mention that 
to define them one should introduce into the game certain configuration  
spaces $D^{\nu}(3)$ whose (more or less obvious) definition we leave 
to the reader. 

The unit of this tensor structure is the sheaf $\One=\CL(0)$ 
(cf. \ref{stand}). 
  
Equipped with these complementary structures, the category $\tFS$ becomes 
a {\em braided tensor category}.

\section{Vanishing cycles}

{\em GENERAL GEOMETRY} 

\subsection{} 
\label{facets} Let us fix a finite set $J$, and consider the space 
$D^J$. Inside this space, let us consider the subspaces $D^J_{\BR}=
D^J\cap\BR^J$ and $D^{J+}=D^J\cap\BR^J_{\geq 0}$. 

Let $\CH$ be the set ({\em arrangement}) 
of all real hyperplanes in $D^J_{\BR}$ of the form 
$H_j:\ t_j=0$ or $H_{j'j''}:\ t_{j'}=t_{j''}$. An {\em edge} $L$ of the 
arrangement $\CH$ is 
a subspace of $D^J_{\BR}$ which is a non-empty intersection $\bigcap H$  
of some hyperplanes from $\CH$. We denote by $L^{\circ}$ the complement 
$L-\bigcup L'$, the union over all edges $L'\subset L$ of smaller 
dimension.   
A {\em facet} of $\CH$ is a connected 
component $F$ of some $L^{\circ}$. We call a facet {\em positive} if 
it lies entirely inside $D^{J+}$.

For example, we have a unique smallest facet $O$ --- the origin. For each 
$j\in J$, we have a positive one-dimensional facet $F_j$ given 
by the equations $t_{j'}=0\ (j'\neq j);\ t_j\geq 0$.  

Let us choose a point $w_F$ on each positive facet $F$. We call a 
{\em flag} a sequence of embedded positive facets $\bF:\ F_0\subset\ldots F_p$; 
we say that $\bF$ {\em starts} from $F_0$. To such a flag we assign 
the simplex $\Delta_{\bF}$ --- the convex hull of the points 
$w_{F_0},\ldots, w_{F_p}$. 

To each positive facet $F$ we assign the following two spaces: 
$D_F=\bigcup\Delta_{\bF}$, the union over all flags $\bF$ starting from $F$, 
and $S_F=\bigcup\Delta_{\bF'}$, the union over all flags $\bF'$ starting 
from a facet which properly contains $F$. Obviously, $S_F\subset D_F$. 

\subsection{} Given a complex $\CK\in\CD(D^J;\CS)$ and a positive 
facet $F$, we introduce a complex of vector spaces $\Phi_F(\CK)$ by 
$$
\Phi_F(\CK)=R\Gamma(D_F,S_F;\CK)[-\dim F]. 
$$
This is a well defined object of the bounded derived category $\CD(*)$ of 
finite dimensional vector spaces,   
not depending on the choice of points $w_F$. It is called the {\em complex 
of vanishing cycles of $\CK$ across $F$}. 

\subsection{Theorem} 
\label{verdier} {\em We have canonically $\Phi_F(D\CK)=D\Phi_F(\CK)$ 
where $D$ denotes the Verdier duality in the corresponding derived 
categories.} $\Box$ 

\subsection{Theorem} {\em If $\CM\in\CM(D^J;\CS)$ then  
$H^i(\Phi_F(\CM))=0$ for $i\neq 0$. Thus, $\Phi_F$ induces an exact functor 
$$
\Phi_F: \CM(D^J;\CS)\lra\Vect
$$
to the category $\Vect$ of finite dimensional vector spaces.} $\Box$ 

\subsection{} 
\label{dir sum} Given a positive facet $E$  
and $\CK\in\CD(D^J,\CS)$,  
we have $S_E=\bigcup_{F\in\CF^1(E)} D_{F}$, the union over the set 
$\CF^1(E)$   
of all positive facets $F\supset E$ with $\dim F=\dim F+1$, and 
$$
R\Gamma(S_E,\bigcup_{F\in\CF^1(E)}\ S_F;\CK)=\oplus_{F\in\CF^1(E)}\ 
R\Gamma(D_F,S_F;\CK).
$$

\subsection{} For two positive facets $E$ and $F\in\CF^1(E)$, and 
$\CK(D^J;\CS)$, let us define the natural map
$$
u=u^F_E(\CK): \Phi_F(\CK)\lra\Phi_E(\CK)
$$  
called {\em canonical}, as the composition
$$
R\Gamma(D_F,S_F;\CK)[-p]\lra R\Gamma(S_E,\bigcup_{F'\in\CF^1(E)} S_{F'};
\CK)[-p]\lra R\Gamma(S_E;\CK)[-p]\lra 
$$
$$
\lra R\Gamma(D_E,S_E)[-p+1]
$$
where $p=\dim F$, the first arrow being induced by the equality in 
\ref{dir sum}, the last one being the coboundary map. 

Define the natural map 
$$
v=v^E_F(\CK): \Phi_E(\CK)\lra\Phi_F(\CK)
$$
called {\em variation}, as the map dual to the composition
$$
D\Phi_F(\CK)=\Phi_F(D\CK)\overset{u(D\CK)}{\lra}\Phi_E(D\CK)=D\Phi_E(\CK).
$$

{\em BACK TO FACTORIZABLE SHEAVES} 

\subsection{} Let $\nu\in Y$. We are going to give two equivalent 
definitions of an exact functor, called {\em vanishing cycles at the origin}
$$
\Phi: \CM(D^{\nu};\CS)\lra\Vect.
$$
{\em First definition.} Let $f: D^{\nu}\lra D$ be the function 
$f((t_j))=\sum t_j$. For an object $\CK\in\CD(D^{\nu};\CS)$, 
the Deligne's complex of vanishing cycles $\Phi_f(\CK)$ 
(cf. \cite{d3}) is concentrated at the origin of the hypersurface $f^{-1}(0)$. 
It is $t$-exact with respect to the middle $t$-structure. We set 
by definition, $\Phi(\CM)= H^0(\Phi_f(\CM))\ (\CM\in\CM(D^{\nu};\CS))$. 

{\em Second definition.} Choose an unfolding of $\nu$, $\pi: J\lra I$. 
Let us consider the canonical projection $\pi: D^J\lra D^{\nu}$. For 
$\CK\in\CD(D^{\nu};\CS)$, the complex $\pi^*\CK$  
is well defined as an element of the $\Sigma_{\pi}$-equivariant 
derived category, hence $\Phi_O(\pi^*(\CK))$ is a well defined object 
of the $\Sigma_{\pi}$-equivariant derived category of vector spaces 
($O$ being the origin facet in $D^J$). 
Therefore, the complex  of $\Sigma_{\pi}$-invariants $\Phi_O(\pi^*\CK)
^{\Sigma_{\pi}}$ is a well defined object of $\CD(*)$. If 
$\CK\in\CM(D^{\nu};\CS)$ then all the cohomology of 
$\Phi_O(\pi^*\CK)^{\Sigma_{\pi}}$ in non-zero degree vanishes.     
We set 
$$
\Phi(\CM)=H^0(\Phi_O(\pi^*\CM)^{\Sigma_{\pi}})\ 
$$
$(\CM\in\CM(D^{\nu};\CS))$. The equivalence of the two definitions follows 
without difficulty from the proper base change theorem. In computations 
the second definition is used.\footnote{its independence of the choice 
of an unfolding follows from its equivalence to the first definition.}

\subsection{} Let $\CM$ be a factorizable sheaf supported at $c\in X/Y$, 
$\lambda=\lambda(\CM)$.  
For $\nu\in Y^+$, define a vector space $\Phi(\CM)_{\lambda-\nu}$ by 
$$
\Phi(\CM)_{\lambda-\nu}=\Phi(\CM^{\nu}).
$$
If $\mu\in X,\ \mu\not\leq\lambda$, set $\Phi(\CM)_{\mu}=0$. One sees easily 
that this way we get an exact functor $\Phi$ from $\tFS_c$ to the category 
of $X$-graded vector spaces with finite dimensional components. 
We extend it to the whole category $\tFS$ by additivity. 

\subsection{} Let $\CM\in\tFS_c$, $\lambda=\lambda(\CM)$, 
$\nu=\sum\nu_ii\in Y^+$. Let $i\in I$ be such that $\nu_i>0$. Pick an 
unfolding of $\nu$, $\pi: J\lra I$. For each $j\in\pi^{-1}(i)$, 
the restriction of $\pi$, $\pi_j: J-\{j\}\lra I$, is an unfolding of 
$\nu-i$. 

For each $j\in\pi^{-1}(i)$, we have canonical and variation morphisms
$$
u_j: \Phi_{F_j}(\pi^*\CM^{\nu})\rlh\Phi_O(\pi^*\CM^{\nu}): v_j
$$
(the facet $F_j$ has been defined in \ref{facets}). 
Taking their sum over $\pi^{-1}(i)$, we get the maps 
$$
\sum u_j:\ \oplus_{j\in\pi^{-1}(i)}\Phi_{F_j}(\pi^*\CM^{\nu})\rlh 
\Phi_O(\pi^*\CM):\sum v_j
$$
Note that the group $\Sigma_{\pi}$ acts on both sides and the maps 
respect this action. After passing to $\Sigma_{\pi}$-invariants, we get 
the maps

(a) $\Phi_{F_j}(\pi^*\CM^{\nu})^{\Sigma_{\pi_j}}=(\oplus_{j'\in\pi^{-1}(i)}
\Phi_{F_{j'}}(\pi^*\CM^{\nu}))^{\Sigma_{\pi}}\rlh 
\Phi_O(\pi^*\CM^{\nu})^{\Sigma_{\pi}}=\Phi(\CM)_{\nu}$. 

Here $j\in\pi^{-1}(i)$ is an arbitrary element. Let us consider the space
$$
F_j^{\perp}=\{(t_{j'})\in D^J|t_j=r, t_{j'}\in D(r')\mbox{ for all }
j'\neq j\} 
$$
where $r, r'$ are some fixed real numbers such that $0<r'<r<1$. 
The space $F^{\perp}_j$ is transversal to $F_j$ and may be identified 
with $D^{J-\{j\}}$. The factorization isomorphism induces the isomorphism 
$$
\pi^*\CM^{\nu}|_{F^{\perp}_j}\cong\pi^*_j\CM^{\nu-i}\otimes(\CI_{\lambda-\nu+i}
^i)_{\{r\}}=\pi^*\CM^{\nu-i}
$$
which in turn induces the isomorphism
$$
\Phi_{F_j}(\pi^*\CM^{\nu})\cong\Phi_O(\pi^*_j\CM^{\nu-i}).
$$
which is $\Sigma_{\pi_j}$-equivariant. Taking $\Sigma_{\pi_j}$-invariants 
and composing with the maps (a), we get the maps 

(b) $\epsilon_i:\Phi(\CM)_{\nu-i}=\Phi_O(\pi^*_j\CM^{\nu-i})^{\Sigma_j}
\rlh \Phi(\CM)_{\nu}: \theta_i$ 

which do not depend on the choice of $j\in\pi^{-1}(i)$. 

\subsection{} For an arbitrary $\CM\in\tFS$,  
let us define the $X$-graded vector space $\Phi(\CM)$ as  
$\Phi(\CM)=\oplus_{\lambda\in X}\Phi(\CM)_{\lambda}$.     

\subsection{Theorem} {\em The operators $\epsilon_i, \theta_i\ (i\in I)$ 
acting on the $X$-graded vector space $\Phi(\CM)$ satisfy the relations 
~\ref{def c} {\em (a), (b)}.} $\Box$ 

\subsection{} A factorizable sheaf $\CM$ is finite 
iff the space $\Phi(\CM)$ is finite dimensional.   
The previous theorem says that $\Phi$ defines an exact functor 
$$
\Phi: \FS\lra\CC.
$$
One proves that $\Phi$ is a tensor functor.  

\subsection{Example.} For every $\lambda\in X$, the factorizable sheaf 
$\CL(\lambda)$ (cf. \ref{stand}) is finite. It is an irreducible object 
of $\FS$, and every irreducible object in $\FS$ is isomorphic to some 
$\CL(\lambda)$. We have $\Phi(\CL(\lambda))=L(\lambda)$.   

The next theorem is the main result of the present work. 

\subsection{Theorem} {\em The functor $\Phi$ is an equivalence of  
braided tensor categories.} $\Box$ 

\subsection{Remark} As a consequence, the category $\FS$ is rigid. We do not  
know a geometric construction of the rigidity; it would be very interesting 
to find one.

\newpage
\begin{center}
{\bf Part II. Global (genus $0$).}
\end{center}

\section{Cohesive local systems} 

\subsection{} From now on until the end of the paper we make the assumptions 
of \ref{k big}. In the operadic notations below we partially follow 
\cite{bd}. 

\subsection{Operad $\CD$} For a nonempty finite 
set $J$, let $D(J)$ denote the space 
whose points are $J$-tuples $\{x_j,\tau_j\}\ (j\in J)$ where $x_j\in D$ and 
$\tau_j$ is a non-zero tangent vector at $x_j$, such that all points $x_j$ 
are distinct. 

Let $\tD(J)$ be the space whose points are $J$-tuples 
$\{\phi_j\}$ of holomorphic maps $\phi_j: D\lra D\ (j\in J)$, 
each $\phi_j$ having the form $\phi_j(z)=x_j+\tau_jz\ (x_j\in D, \tau_j\in
\BC^*)$, such that $\phi_j(D)\cap\phi_{j'}(D)=\emptyset$ for 
$j\neq j'$. We shall identify the $j$-tuple $\{\phi_j\}$ with the 
$J$-tuple $\{x_j,\tau_j\}$, and consider $\tau_j$ as a non-zero tangent 
vector from $T_{x_j}D$, thus identifying $T_{x_j}D$ with $\BC$ using 
the local coordinate $z-x_j$. So, $\tau_j$ is the image under $\phi_j$ of 
the unit tangent vector at $0$.   
We have an obvious map $p(J):\tD(J)\lra D(J)$ 
which is a homotopy equivalence.  

If $\rho: K\lra J$ is an epimorphic map of finite sets, the composition defines 
a holomorphic map 
$$
m(\rho):\ \prod_J \tD(K_j)\times \tD(J)\lra \tD(K)
$$
where $K_j:=\rho^{-1}(j)$. 
If $L\overset{\sigma}{\lra}K\overset{\rho}{\lra}J$ are two epimorphisms of finite 
sets, the square
$$\begin{array}{ccc}
\prod_K \tD(L_k)\times\prod_J \tD(K_j)\times \tD(J)&\overset{m(\sigma)}{\lra}&
\prod_K \tD(L_k)\times\tD(K)\\
\prod m(\sigma_j)\downarrow&\ &\downarrow m(\sigma)\\
\prod_J \tD(L_j)\times\tD(J)&\overset{m(\rho\sigma)}{\lra}&\tD(L)
\end{array}$$
commutes. Here $\sigma_j: L_j\lra K_j$ are induced by $\sigma$. 

Let $*$ denote the one element set. 
The space $\tD(*)$ has a marked point, also to be denoted by $*$, corresponding 
to the identity map $\phi: D\lra D$. 

If $\rho: J'\iso J$ is an isomorphism, it induces in the obvious way 
an isomorphism $\rho^*:\tD(J)\iso\tD(J')$ (resp. $D(J)\iso D(J')$). The first 
map coincides with $m(\rho)$ restricted to $(\prod_J *)\times\tD(J)$. 
In particular, 
for each $J$, the group $\Sigma_J$ of automorphisms of the set $J$, acts 
on the spaces $\tD(J), D(J)$.   

The map $m(J\lra *)$ restricted to $*\times\tD(J)$ is the identity of $\tD(J)$. 

We will denote the collection of the spaces and maps $\{\tD(J), m(\rho)\}$ 
by $\CD$, and call it the {\em operad of disks with tangent vectors}. 

\subsection{Coloured local systems over $\CD$} If $\rho: K\lra J$ 
is an epimorphism of finite sets and $\pi: K\lra X$ is a map of sets, 
we define the map $\rho_*\pi: J\lra X$ by $\rho_*\pi(j)=\sum_{K_j}\pi(k)$. 
For $j\in J$, we denote $K_j:=\rho^{-1}(j)$ as above, and $\pi_j:\ K_j\lra X$ 
will denote the restriction of $\pi$.  

Let us call an {\em $X$-coloured local system} $\CJ$ over $\CD$ 
a collection of local systems $\CJ(\pi)$ over the spaces $\tD(J)$ given 
for every map $\pi:J\lra X$, $J$ being a non-empty finite set, together 
with {\em factorization isomorphisms}
$$
\phi(\rho):\ m(\rho)^*\CJ(\pi)\iso\Boxtimes_J\CJ(\pi_j)\boxtimes\CJ(\rho_*\pi)
$$
given for every epimorphism $\rho: K\lra J$ and $\pi: K\lra X$, which 
satisfy the properties (a), (b) below.  

(a) {\em Associativity}. Given a map $\pi: L\lra X$ and a pair of epimorphisms 
$L\overset{\sigma}{\lra}K\overset{\rho}{\lra}J$, the square below commutes.
$$\begin{array}{ccc}
m(\rho)^*m(\sigma)^*\CJ(\pi)&\overset{\phi(\sigma)}{\lra}&\Boxtimes_K\CJ(\pi_k)
\boxtimes m(\rho)^*\CJ(\sigma_*\pi)\\
\phi(\rho\sigma)\downarrow&\ &\downarrow\phi(\rho)\\
\Boxtimes_J m(\sigma_j)^*\CJ(\pi_j)\boxtimes\CJ(\rho_*\sigma_*\pi)&\overset
{\boxtimes\phi(\sigma_j)}{\lra}&\Boxtimes_K\CJ(\pi_k)\boxtimes\Boxtimes_J
\CJ((\sigma_*\pi)_j)\boxtimes\CJ(\rho_*\sigma_*\pi)
\end{array}$$
Note that $\pi_{j*}\sigma_j=(\sigma_*\pi)_j$. 

For $\mu\in X$, let $\pi_{\mu}: *\lra X$ be defined by $\pi_{\mu}(*)=\mu$. 
The isomorphisms $\phi(\id_*)$ restricted to the marked points 
in $D(*)$, give the isomorphisms $\CJ(\pi_{\mu})_*\iso k$ (and imply that 
the local systems $\CJ(\pi_{\mu})$ are one-dimensional). 

(b) For any $\pi: J\lra X$, the map $\phi(\id_J)$ restricted to 
$(\prod_J *)\times \tD(J)$, equals $\id_{\CJ(\pi)}$. 

The map $\phi(J\lra *)$ restricted to $*\times\tD(J)$, equals $\id_{\CJ(\pi)}$. 

\subsection{} The definition above implies that the local systems $\CJ(\pi)$ 
are functorial with respect to isomorphisms. In particular, 
the action of the group $\Sigma_{\pi}$ on $\tD(J)$ lifts to $\CJ(\pi)$. 

\subsection{Standard local system over $\CD$} 
\label{stand d} Let us define the "standard" local systems 
$\CJ(\pi)$ by a version of the construction \ref{stand loc}.   

We will use 
the notations of {\em loc. cit.} for the marked points 
and paths. We can identify $D(J)=D^{J\circ}\times(\BC^*)^J$ where 
we have identified all tangent spaces $T_xD\ (x\in D)$ with $\BC$ using 
the local coordinate $z-x$. 

We set $\CJ(\pi)_{x(\sigma)}=k$, the monodromies being 
$T_{\gamma(j',j'')^{\pm}}=\zeta^{\mp\pi(j')\cdot\pi(j'')}$, and the monodromy 
$T_j$  
corresponding to the counterclockwise circle of a tangent vector $\tau_j$ 
is $\zeta^{-2n\pi(j)}$. The factorization isomorphisms 
are defined by the same condition as in {\em loc. cit}. 

This defines 
the {\em standard local system $\CJ$} over $\CD$. Below, the notation $\CJ$ 
will be reserved for this local system.

\subsection{} 
\label{right} Set $P=\BP^1(\BC)$. We pick a point $\infty\in P$ and 
choose a global coordinate $z:\BA^1(\BC)=P-\{\infty\}\iso\BC$. 
This gives local coordinates: $z-x$ at $x\in\BA^1(\BC)$ and 
$1/z$ at $\infty$. 

For a non-empty finite 
set $J$, let $P(J)$ denote the space of $J$-tuples $\{ x_j,\tau_j\}$ where 
$x_j$ are distinct points on $P$, and $\tau_j$ is a non-zero tangent vector 
at $x_j$. Let $\tP(J)$ denote the space whose points are $J$-tuples of 
holomorphic embeddings $\phi_j: D\lra P$ with non-intersecting images 
such that each $\phi_j$ is a restriction of an algebraic morphism $P\lra P$. 
We have the $1$-jet projections $\tP(J)\lra P(J)$. We will use the notation 
$P(n)$ for $P(\{1,\ldots,n\})$, etc.   

An epimorphism $\rho: K\lra J$ induces the maps 
$$
m_P(\rho):\ \prod_J \tD(K_j)\times \tP(J)\lra\tP(K)
$$
and
$$
\bar{m}_P(\rho):\ \prod_J D(K_j)\times\tP(J)\lra P(K).
$$
For a pair of epimorphisms $L\overset{\sigma}{\lra}K\overset{\rho}{\lra}J$, 
the square
$$\begin{array}{ccc}
\prod_K \tD(L_k)\times\prod_J \tD(K_j)\times \tP(J)&\overset{m_P(\sigma)}{\lra}&
\prod_K \tD(L_k)\times\tP(K)\\
\prod m(\sigma_j)\downarrow&\ &\downarrow m_P(\sigma)\\
\prod_J \tD(L_j)\times\tP(J)&\overset{m_P(\rho\sigma)}{\lra}&\tP(L)
\end{array}$$
commutes. 

If $\rho: J'\iso J$ is an isomorphism, it induces in the obvious way 
an isomorphism $\rho^*:\tP(J)\iso\tP(J')$ (resp. $P(J)\iso P(J')$). This 
last map 
coincides with $m_P(\rho)$ (resp. $\bar{m}_P(\rho)$)  
restricted to $(\prod_J *)\times\tP(J)$ (resp. $(\prod_J *)\times P(J)$). 
 
The collection of the spaces and maps $\{\tP(J), m_P(\rho)\}$ form an 
object $\tP$ called a {\em right module} over the operad $\CD$. 

\subsection{} 
\label{cls def} Fix an element $\mu\in X$. Let us say that a map $\pi: J\lra X$ 
has {\em level $\mu$} if $\sum_J\pi(j)=\mu$.  

A {\em cohesive local 
system $\CH$ of level $\mu$} over $P$  
is a collection of local systems $\CH(\pi)$ over 
$\tP(J)$ given for every $\pi: J\lra X$ of level $\mu$, 
together with {\em factorization isomorphisms}
$$
\phi_P(\rho): m_P(\rho)^*\CH(\pi)\iso\Boxtimes_J\CJ(\pi_j)\boxtimes
\CH(\rho_*\pi)
$$
given for every epimorphism $\rho: K\lra J$ and $\pi: K\lra X$ of level 
$\mu$, which satisfy the properties (a), (b) below. 

(a) {\em Associativity}. Given a map $\pi: L\lra X$ and a pair of epimorphisms 
$L\overset{\sigma}{\lra}K\overset{\rho}{\lra}J$, the square below commutes.
$$\begin{array}{ccc}
m_P(\rho)^*m_P(\sigma)^*\CH(\pi)&\overset{\phi_P(\sigma)}{\lra}&
\Boxtimes_K\CJ(\pi_k)
\boxtimes m_P(\rho)^*\CH(\sigma_*\pi)\\
\phi_P(\rho\sigma)\downarrow&\ &\downarrow\phi_P(\rho)\\
\Boxtimes_J m(\sigma_j)^*\CJ(\pi_j)\boxtimes\CH(\rho_*\sigma_*\pi)&\overset
{\boxtimes\phi(\sigma_j)}{\lra}&\Boxtimes_K\CJ(\pi_k)\boxtimes\Boxtimes_J
\CJ((\sigma_*\pi)_j)\boxtimes\CH(\rho_*\sigma_*\pi)
\end{array}$$

(b) For any $\pi: J\lra X$, the map $\phi_P(\id_J)$ restricted to 
$(\prod_J *)\times \tP(J)$, equals $\id_{\CH(\pi)}$. 

\subsection{} The definition above implies that the local systems $\CH(\pi)$ 
are functorial with respect to isomorphisms. In particular, 
the action of the group $\Sigma_{\pi}$ on $\tP(J)$, lifts to $\CH(\pi)$. 

\subsection{Theorem} {\em For each $\mu\in X$ such that 
$\mu\equiv\ 2\rho_{\ell}\ \modul Y_{\ell}$, there exists a unique 
up to an isomorphism one-dimensional cohesive local system $\CH=\CH^{(\mu)}$ 
of level $\mu$ over $P$}. $\Box$ 

The element $\rho_{\ell}\in X$ is defined in \ref{l} and the lattice 
$Y_{\ell}$ in \ref{lattice}.   

>From now on, let us fix such a local system $\CH^{(\mu)}$ for each $\mu$ as 
in the theorem. 

\subsection{} Note that the obvious maps $p(J): \tP(J)\lra P(J)$ are homotopy 
equivalences. Therefore the local systems $\CH(\pi)\ (\pi: J\lra X)$ 
descend to the unique local systems over $P(J)$, to be denoted by the same 
letter.

\section{Gluing} 

\subsection{} 
\label{diag} Let us fix a finite set $K$. For $\nu\in Y^+$, pick 
an unfolding of $\nu$, $\pi: J\lra I$. Consider the space $P^{\nu}$; 
its points are formal linear combinations $\sum_J \pi(j) x_j\ 
(x_j\in P)$. We define the space $P^{\nu}(K)=P(K)\times P^{\nu}$; its points 
are tuples $\{y_k,\tau_k,\sum\pi(j)x_j\}\ (k\in K, y_k\in P, \tau_k\neq 0$ in  
$T_{y_k}P)$, all $y_k$ being distinct. Let $P^{\nu}(K)^{\circ}$ 
(resp. $P^{\nu}(K)^{\bullet}$)  
be the subspace whose points are tuples as above with all $x_j$ distinct 
from $y_k$ and pairwise distinct (resp. all $x_j$ distinct from $y_k$). 
We will use the notation $P^{\nu}(n)$ for $P^{\nu}(\{1,\ldots,n\})$, etc.  

Let $\nu'\in Y^+$ and $\vnu=\{\nu_k\}\in (Y^+)^K$ be such that 
$\sum_K\nu_k+\nu'=\nu$. Define the space $P^{\vnu;\nu'}\subset\tP(K)\times 
P^{\nu}$ consisting of tuples $\{\phi_k,\sum\pi(j)x_j\}$ such that 
for each $k$, $\nu_k$ of the points $x_j$ lie inside $\phi_k(D)$ and 
$\nu'$ of them lie outside all closures of these disks. 

Let $\vzero\in(Y^+)^K$ be the zero $K$-tuple. Define the space 
$\tP^{\nu}(K):=P^{\vzero;\nu}$.  
We have obvious maps 

(a) $\prod_K D^{\nu_k}\times\tP^{\nu'}(K)\overset{p(\vnu;\nu')}{\lla}
P^{\vnu;\nu'}\overset{m(\vnu;\nu')}{\lra}P^{\nu}(K)$. 

Let $\vnu^1, \vnu^2\in(Y^+)^K$ be such that $\vnu^1+\vnu^2=\vnu$. Define the 
space 

$P^{\vnu^1;\vnu^2;\nu}\subset\tP(K)\times\BR^K_{>0}\times\tP(K)
\times P^{\nu}$ consisting of all tuples $\{\phi_k; r_k; \phi'_k; 
\sum\pi(j)x_j\}$ such that $r_k<1; \phi_k(z)=\phi'_k(r_kz)$, 
and $\nu^1_k$ (resp. $\nu^2_k$, $\nu'$) from the points $x_j$ lie inside 
$\phi_k(D)$ (resp. inside the annulus $\phi'_k(D)-\overline{\phi_k(D)}$, 
inside $P-\bigcup\overline{\phi'_k(D)}$). We set $\tP^{\vnu;\nu'}:=
P^{\vzero;\vnu;\nu'}$, cf. \ref{spaces}. 

We have a commutative romb (cf. \ref{rhomb}):   

(b) \begin{center}
  \begin{picture}(14,8)
    \put(7,8){\makebox(0,0){$P^{\nu}(K)$}}
    \put(5,6){\makebox(0,0){$P^{\vnu^1;\nu^2+\nu'}$}}
    \put(9,6){\makebox(0,0){$P^{\vnu^1+\vnu^2;\nu'}$}}
    \put(3,4){\makebox(0,0){$\prod D^{\nu^1_k}\times\tP^{\nu^2+\nu'}(K)$}}
    \put(7,4){\makebox(0,0){$P^{\vnu^1;\vnu^2;\nu'}$}}
    \put(11,4){\makebox(0,0){$\prod D^{\nu^1_k+\nu^2_k}\times\tP^{\nu'}(K)$}}
    \put(5,2){\makebox(0,0){$\prod D^{\nu^1_k}\times\tP^{\vnu^2;\nu'}$}}
    \put(9,2){\makebox(0,0){$\prod D^{\nu^1_k,\nu^2_k}\times\tP^{\nu'}(K)$}}
    \put(7,0){\makebox(0,0){$\prod D^{\nu^1_k}\times\prod A^{\nu^2_k}
              \times\tP^{\nu'}(K)$}}

    \put(5.5,6.5){\vector(1,1){1}}
    \put(8.5,6.5){\vector(-1,1){1}}
    \put(4.5,5.5){\vector(-1,-1){1}}
    \put(6.5,4.5){\vector(-1,1){1}}
    \put(7.5,4.5){\vector(1,1){1}}
    \put(9.5,5.5){\vector(1,-1){1}} 
    \put(4.5,2.5){\vector(-1,1){1}}
    \put(6.5,3.5){\vector(-1,-1){1}}
    \put(7.5,3.5){\vector(1,-1){1}}
    \put(9.5,2.5){\vector(1,1){1}}
    \put(5.5,1.5){\vector(1,-1){1}}
    \put(8.5,1.5){\vector(-1,-1){1}}

    \put(5.5,7){\makebox(0,0){$m$}}
    \put(8.5,7){\makebox(0,0){$m$}}
    \put(3.5,5){\makebox(0,0){$p$}}
    \put(5.5,5){\makebox(0,0){$m$}}
    \put(8.5,5){\makebox(0,0){$m$}}
    \put(10.5,5){\makebox(0,0){$p$}}
    \put(3.5,3){\makebox(0,0){$m$}}
    \put(5.5,3){\makebox(0,0){$p$}}
    \put(8.5,3){\makebox(0,0){$p$}}
    \put(10.5,3){\makebox(0,0){$m$}}
    \put(5.5,1){\makebox(0,0){$p$}}
    \put(8.5,1){\makebox(0,0){$p$}}
    
  \end{picture} 
\end{center}

Here $\nu^2:=\sum\nu^2_k$.

\subsection{} Let $P(K;J)$ be the space consisting of tuples $\{y_k;\tau_k;
x_j\}\ (k\in K, j\in J, y_k, x_j\in P, \tau_k\neq 0$ in $T_{y_k}P$), where 
all points $x_k, y_j$ are distinct. We have 
$P^{\nu}(K)^{\circ}=P(K;J)/\Sigma_{\pi}$.  
We have an obvious projection $P(K\coprod J)\lra P(K;J)$.

\subsection{} 
\label{sign p} Let $\{\CM_k\}$ be a $K$-tuple of factorizable sheaves 
supported at some cosets in $X/Y$; let $\mu_k=\lambda(\CM_k)$.   

Let $\tpi: K\coprod J\lra X$ be a map defined by $\tpi(k)=\mu_k$, 
$\tpi(j)=-\pi(j)\in I\hra X$. The local system $\CH(\tpi)$ over $P(K\coprod J)$ 
descends to $P(K;J)$ since $\zeta^{2n(-i)}=1$ for all $i\in I$;  
this one in turn descends to the unique local system 
$\tCH_{\vmu}^{\nu}$ over $P^{\nu}(K)^{\circ}$, 
due to $\Sigma_{\pi}$-equivariance. Let us define the local system 
$\CH_{\vmu}^{\nu}:=\tCH_{\vmu}^{\nu}\otimes\Sign^{\nu}$.  
Here $\Sign^{\nu}$ denotes the inverse image of the sign local 
system on $P^{\nu\circ}$ (defined in the same manner as for the disk, cf. 
\ref{sign}) under the forgetful map $P^{\nu}(K)^{\circ}\lra P^{\nu\circ}$. 
    
Let $\CH^{\nu\bullet}_{\vmu}$ 
be the perverse sheaf over $P^{\nu}(K)^{\bullet}$ which is  
the middle extension of $\CH^{\nu}_{\vmu}[\dim P^{\nu}(K)]$.   
Let us denote by the same letter the inverse image of this perverse sheaf 
on the space $\tP^{\nu}(K)$ with respect to the evident projection 
$\tP^{\nu}(K)\lra P^{\nu}(K)^{\bullet}$. 

\subsection{} 
\label{adm} Let us call an element $\nu\in Y^+$ {\em admissible} 
(for a $K$-tuple $\{\mu_k\}$)  
if $\sum\mu_k-\nu\equiv\ 2\rho_{\ell}\ \modul Y_{\ell}$, see \ref{lattice}.  

\subsection{Theorem - definition} 
\label{glu thm} {\em For each admissible $\nu$, there exists 
a unique, up to a unique isomorphism, perverse sheaf, denoted by\  
$\Boxtimes^{(\nu)}_K\ \CM_k$, over $P^{\nu}(K)$, equipped with 
isomorphisms 
$$
\psi(\vnu;\nu'): m^*\Boxtimes^{(\nu)}_K\ \CM_k
\iso p^*(\Boxtimes_K\ \CM_k^{\nu_k}
\boxtimes\CH^{\nu'\bullet}_{\vmu-\vnu})
$$
given for every diagram \ref{diag} (a) such that for each rhomb 
\ref{diag} (b) the cocycle condition
$$
\phi(\vnu^2;\nu')\circ\psi(\vnu^1;\nu^2+\nu')=
(\Boxtimes_K\ \psi^{\CM_k}(\nu^1_k,\nu^2_k))\circ\psi(\vnu^1+\vnu^2;\nu')
$$
holds.} $\Box$

\subsection{} The sheaf $\Boxtimes^{(\nu)}_K\ \CM_k$ defines for each $K$-tuple 
of $\vy=\{y_k,\tau_k\}$ of points of $P$ with non-zero tangent vectors, 
the sheaf $\Boxtimes^{(\nu)}_{\vy}\ \CM_k$ over $P^{\nu}$, to be called the 
{\em gluing of the factorizable sheaves $\CM_k$ into the points $(y_k,\tau_k)$}. 

\subsection{Example} 
\label{middle} The sheaf $\Boxtimes^{(\nu)}_K\ \CL(\mu_k)$ is equal 
to the middle extension of the sheaf $\CH_{\vmu}^{\nu\bullet}$.

\section{Semiinfinite cohomology} 

In this section we review the theory of semiinfinite cohomology 
in the category $\CC$, due to S. M. Arkhipov, cf. \cite{arkh}.  

\subsection{} Let $\CC_r$ be a category whose objects are {\em right} 
$\fu$-modules $N$, finite dimensional over $k$, with a given $X$-grading 
$N=\oplus_{\lambda\in X}\ N_{\lambda}$ such that 
$xK_{\nu}=\zeta^{-\langle\nu,\lambda\rangle}x$ for any $\nu\in Y, 
\lambda\in X, x\in N_{\lambda}$. All definitions and results  
concerning the 
category $\CC$ given above and below, have the obvious versions 
for the category $\CC_r$. 

For $M\in\CC$, define $M^{\vee}\in\CC_r$ as follows: $(M^{\vee})_{\lambda}=
(M_{-\lambda})^*$ (the dual vector space); the action of the operators 
$\theta_i, \epsilon_i$ being the transpose of their action on $M$. This 
way we get an equivalence $^{\vee}:\CC^{\opp}\iso\CC_r$.         

\subsection{} Let us call an object $M\in\CC$ $\fu^-$- (resp. $\fu^+$-) 
{\em good} if it admits a filtration whose successive quotients 
have the form $\ind_{\fu^{\geq 0}}^{\fu}(M')$ (resp. 
$\ind_{\fu^{\leq 0}}^{\fu}(M'')$) for some $M'\in\CC^{\geq 0}$ 
(resp. $M''\in\CC^{\leq 0}$) (cf. \ref{ind}). These classes of objects  
are stable with respect to the tensor multiplication by an arbitrary object 
of $\CC$.  

If $M$ is $\fu^-$- (resp. $\fu^+$-) good then $M^*$ is $\fu^-$- 
(resp. $\fu^+$-) good. 

If $M$ is $\fu^-$-good and $M'$ is $\fu^+$-good then $M\otimes M'$ is a 
projective object in $\CC$. 

\subsection{} Let us say that a complex $M^{\bullet}$ in $\CC$ is 
{\em concave} (resp. {\em convex}) if 

(a) there exists $\mu\in X$ such that all nonzero components 
$M^{\bullet}_{\lambda}$ have the weight $\lambda\geq\mu$ (resp. 
$\lambda\leq\mu$);  

(b) for any $\lambda\in X$, the complex $M^{\bullet}_{\lambda}$ is finite. 

\subsection{} For an object $M\in\CC$, we will call a {\em left} 
(resp. {\em right}) {\em $\fu^{\pm}$-good resolution} of $M$  an exact 
complex 
$$
\ldots\lra P^{-1}\lra P^0\lra M\lra 0
$$
(resp.
$$ 
0\lra M\lra R^0\lra R^1\lra\ldots)
$$
such that all $P^i$ (resp. $R^i$) are $\fu^{\pm}$-good. 

\subsection{Lemma.} {\em Each object $M\in\CC$ admits a convex $\fu^-$-good 
left resolution, a concave $\fu^+$-good left resolution, a concave 
$\fu^-$-good right resolution and a convex $\fu^+$-good right resolution.} 
$\Box$ 

\subsection{} 
\label{pairing} For $N\in\CC_r, M\in\CC$, define a vector space $N\otimes_{\CC}M$ 
as the zero weight component of the tensor product $N\otimes_{\fu}M$ 
(which has an obvious $X$-grading). 

For $M, M'\in\CC$, we have an obvious perfect pairing 

(a) $\Hom_{\CC}(M,M')\otimes(M^{\prime\vee}\otimes_{\CC}M)\lra k$. 

\subsection{} $M, M'\in\CC, N\in\CC_r$ and $n\in\BZ$, define the 
{\em semiinifinite $\Ext$ and $\Tor$} spaces
$$
\Ext^{\infh+n}_{\CC}(M,M')=H^n(\Hom_{\CC}(R^{\bullet}_{\searrow}(M),
R^{\bullet}_{\nearrow}(M')))
$$
where $R^{\bullet}_{\searrow}(M)$ (resp. $R^{\bullet}_{\nearrow}(M')$) 
is an arbitrary $\fu^+$-good convex right resolution of $M$ 
(resp. $\fu^-$-good concave right resolution of $M'$),
$$
\Tor^{\CC}_{\infh+n}(N,M)=H^{-n}(P^{\bullet}_{\swarrow}(N)\otimes_{\CC}
R^{\bullet}_{\searrow}(M))
$$
where $P^{\bullet}_{\swarrow}(N)$ is an arbitrary $\fu^-$-good convex left 
resolution of $N$. 

This definition does not depend, up to a unique isomorphism, 
upon the choice of resolutions, and is functorial.
  
These spaces are finite dimensional and are non-zero only for finite 
number of degrees $n$. 

The pairing \ref{pairing} (a) induces perfect pairings
$$
\Ext^{\infh+n}_{\CC}(M,M')\otimes\Tor^{\CC}_{\infh+n}(M^{\prime\vee},M)\lra k\ 
(n\in\BZ). 
$$

\section{Conformal blocks (genus $0$)}

In this section we suppose that $k=\BC$. 

\subsection{} 
\label{trivial} Let $M\in \CC$. We have a canonical embedding of vector 
spaces 
$$
\Hom_{\CC}(\One,M)\hra M
$$
which identifies $\Hom_{\CC}(\One,M)$ with the maximal trivial
subobject of $M$. Here "trivial" means "isomorphic to a sum of a few 
copies of the object $\One$".   
Dually, we have a canonical epimorhism
$$
M\lra\Hom_{\CC}(M,\One)^*
$$
which identifies $\Hom_{\CC}(M,\One)^*$ with the maximal trivial quotient
of $M$. Let us denote by $\langle M\rangle$ the image
of the composition
$$
\Hom_{\CC}(\One,M)\lra M\lra\Hom_{\CC}(M,\One)^*
$$
Thus, $\langle M\rangle$ is canonically a subquotient of $M$.

One sees easily that if $N\subset M$ is a trivial direct summand of $M$
which is maximal, i.e. not contained in greater direct summand, then
we have a canonical isomorphism $\langle M\rangle\iso N$. For this reason,
we will call $\langle M\rangle$ {\em the maximal trivial
direct summand} of $M$. 

\subsection{} Let $\gamma_0\in Y$ denote the highest coroot and $\beta_0\in Y$ 
denote the 
coroot dual to the highest root ($\gamma_0=\beta_0$ for a simply laced 
root datum).
 
Let us define the {\em first alcove} $\Delta_{\ell}\subset X$ by 
$$
\Delta_{\ell}=\{\lambda\in X|\ \langle i,\lambda+\rho\rangle>0\mbox{ for all }
i\in I;\ \langle\gamma_0,\lambda+\rho\rangle<\ell\}
$$
if $d$ does not divide $\ell$, i.e. if $\ell_i=\ell$ for all $i\in I$, 
and by
$$
\Delta_{\ell}=\{\lambda\in X|\ \langle i,\lambda+\rho\rangle >0\mbox{ for all }
i\in I;\ \langle\beta_0,\lambda+\rho\rangle <\ell_{\beta_0}\}
$$
otherwise, cf. \cite{ap} 3.19 ($d$ is defined in \ref{root}, 
and $\ell_{\beta_0}$ in \ref{l}). Note that $\ell_{\beta_0}=\ell/d$.   

\subsection{} For $\lambda_1,\ldots,\lambda_n\in\Delta_{\ell}$, define the 
{\em space of conformal blocks} by 
$$
\langle L(\lambda_1),\ldots,L(\lambda_n)\rangle :=
\langle L(\lambda_1)\otimes\ldots\otimes L(\lambda_n)\rangle.
$$
In fact, due to the ribbon structure on $\CC$, the right hand 
side is a {\em local system} over the space $P(n):=
P(\{1,\ldots,n\})$ (cf. \cite{d1}). It is more appropriate to 
consider the previous equality as the definition of the 
{\em local system of conformal blocks} over $P(n)$.

\subsection{Theorem} 
\label{ar} (Arkhipov) {\em For each $\lambda_1,\ldots,\lambda_n
\in\Delta_{\ell}$, the space of conformal blocks  
$\langle L(\lambda_1),\ldots,L(\lambda_n)\rangle$ is naturally a 
subquotient of the space

$\Tor^{\CC}_{\infh+0}({\em\One}_r, L(\lambda_1)\otimes\ldots\otimes 
L(\lambda_n)\otimes L(2\rho_{\ell}))$.

More precisely, due to the ribbon structure on $\CC$, the latter space 
is a stalk 
of a local system over $P(n+1)$, and inverse image of the local system 
$\langle L(\lambda_1),\ldots,L(\lambda_n)\rangle$ under the projection 
onto the first coordinates $P(n+1)\lra P(n)$, is a natural 
subquotient of this local system.} 

Here $\One_r$ is the unit object in $\CC_r$.

Examples, also due to Arkhipov, show that the local systems of conformal blocks 
are in general {\em proper} subquotients of the corresponding  
$\Tor$ local systems.   

This theorem is an immediate consequence of the next lemma, 
which in turn follows from the geometric theorem \ref{integral} below, 
cf. \ref{delta}.  

\subsection{Lemma} 
\label{lemma rho}{\em We have $\Tor^{\CC}_{\infh+n}({\em\One}_r,
L(2\rho_{\ell}))=k$ if $n=0$, and $0$ otherwise.} $\Box$ 

\subsection{} 
\label{charge} Let $\fg$ be the simple Lie algebra (over $k$) associated 
with our Cartan datum; let $\hfg$ be the corresponding affine 
Lie algebra. 

Let $\MS$ denote the category of integrable $\hfg$-modules of central 
charge $\kappa-h$. Here $\kappa$ is a fixed positive integer, $h$ is the 
dual Coxeter number of our Cartan datum. $\MS$ is a semisimple abelian  
category whose irreducible objects are 
$\fL(\lambda),\ \lambda\in\Delta_{\ell}$ where  
$l=2d\kappa$, 
i.e. $\ell=d\kappa$ (we are grateful to Shurik Kirillov who pointed 
out the necessity of the factor $d$ here) and $\fL(\lambda)$ is the 
highest weight module 
with a highest vector $v$ whose "top" part $\fg\cdot v$ is the irreducible 
$\fg$-module of the highest weight $\lambda$. 

According to Conformal field theory, $\MS$ has a natural structure 
of a ribbon category, cf. \cite{ms}, \cite{k}.  

The usual local systems of conformal blocks in the WZW model may be defined as  
$$
\langle\fL(\lambda_1),\ldots,\fL(\lambda_n)\rangle=
\Hom_{\MS}(\One,\fL(\lambda_1)\otimes\ldots\otimes\fL(\lambda_n))
$$
the structure of a local system on the right hand side is due to the ribbon 
structure on $\MS$.  

\subsection{}
\label{ms} Let  $\zeta=\exp(\pi\sqrt{-1}/d\kappa)$.   
We have an exact functor 
$$
\phi:\ \MS\lra\CC 
$$
sending $\fL(\lambda)$ to $L(\lambda)$. This functor identifies 
$\MS$ with a full subcategory of $\CC$. 

The functor $\phi$ does not respect the tensor structures. It admits the 
left and right adjoints, $\phi^{\flat}, \phi^{\sharp}$. For $M\in\CC$, let 
$\langle M\rangle_{\MS}$ denotes the image of the composition 
$$
\phi\circ\phi^{\flat}(M)\lra M\lra\phi\circ\phi^{\sharp}(M).
$$
We have the following comparison theorem.

\subsection{Theorem}{\em We have naturally  
$$
\phi(\fM\otimes\fM')=\langle\phi(\fM)\otimes\phi(\fM')\rangle_{\MS}.
$$}

This follows from the combination of the results of \cite{ap}, \cite{kl}, 
\cite{l3} and \cite{f}. 

\subsection{Corollary.} {\em For any $\lambda_1,\ldots,\lambda_n\in\Delta_{\ell}$, 
the functor $\phi$ induces an isomorphism of local systems 
$$
\langle\fL(\lambda_1),\ldots,\fL(\lambda_n)\rangle=
\langle L(\lambda_1),\ldots,L(\lambda_n)\rangle.\ \Box
$$}

\section{Integration}

We keep the notations of the previous section. 

\subsection{} Let $K$ be a finite set, $m=\card(K)$, $\{\CM_k\}$ a $K$-tuple 
of finite factorizable sheaves, $\CM_k\in\FS_{c_k}$, $\mu_k:=\lambda
(\CM_k)$. Assume that $\nu:=\sum_K\ \mu_k-2\rho_{\ell}$ belongs to $Y^+$.  

Let $\eta: P^{\nu}(K)\lra P(K)$ be the projection. 

\subsection{Theorem} 
\label{integral} {\em We have canonical isomorphisms of local systems 
over $P(K)$ 
$$
R^{a-2m}\eta_*(\Boxtimes^{(\nu)}_K\ \CM_k) =\Tor^{\CC}_{\infh-a}(\em{\One}_r,
\otimes_K\ \Phi(\CM_k))\ (a\in\BZ),
$$
the structure of a local system on the right hand side being induced 
by the ribbon structure on $\CC$.} $\Box$ 

\subsection{Proof of Lemma \ref{lemma rho}} 
\label{delta} Apply the previous theorem 
to the case when the $K$-tuple consists of one sheaf $\CL(2\rho_{\ell})$ and 
$\nu=0$. $\Box$    

\subsection{} From now until the end of the section, $k=\BC$ and 
$\zeta=\exp(\pi\sqrt{-1}/d\kappa)$.    
Let $\lambda_1,\ldots,\lambda_n\in\Delta_{\ell}$. 
Let $\nu=\sum_{m=1}^n\ \lambda_m$; assume that $\nu\in Y^+$. Set   
$\vmu=\{\lambda_1,\ldots,\lambda_n,2\rho_{\ell}\}$. 

Let $\eta$ be the projection   
$P^{\nu}(n+1):=P^{\nu}(\{1,\ldots,n+1\})\lra P(n+1)$ and 
$p: P(n+1)\lra P(n)$ be the projection on the first coordinates. 
 
Let $\CH^{\nu\sharp}_{\vmu}$ denote the middle extension of the sheaf 
$\CH^{\nu\bullet}_{\vmu}$. By the Example \ref{middle}, 
$$
\CH^{\nu\sharp}_{\vmu}= 
\Boxtimes_{1\leq a\leq n}^{(\nu)}\CL(\lambda_a)\boxtimes
\CL(2\rho_{\ell}).
$$   
\subsection{Theorem} 
\label{subq} {\em The local system $p^*\langle L(\lambda_1),\ldots,
L(\lambda_n)\rangle$ is canonically a subquotient of the local 
system 
$$
R^{-2n-2}\eta_*\CH^{\nu\sharp}_{\vmu}.
$$} 

This theorem is an immediate corollary of the previous one and of the Theorem  
\ref{ar}. 

\subsection{Corollary} {\em In the notations of the previous theorem, 
the local system $\langle L(\lambda_1),\ldots,L(\lambda_n)\rangle$ 
is semisimple.} 

{\bf Proof.} The local system 
$p^*\langle L(\lambda_1),\ldots,L(\lambda_n)\rangle$ is a subquotient 
of the geometric local system $R^{-2n-2}\eta_*\CH^{\nu\sharp}_{\vmu}$, and 
hence is semisimple by 
the Beilinson-Bernstein-Deligne-Gabber 
Decomposition theorem, \cite{bbd}, Th\'{e}or\`{e}me 6.2.5. Therefore, 
the local system $\langle L(\lambda_1),\ldots,L(\lambda_n)\rangle$ is 
also semisimple, since the map $p$ induces the surjection on the fundamental 
groups. $\Box$ 

\subsection{} 
\label{hermit} For a sheaf $\CF$, let 
$\bar{\CF}$ denote the sheaf obtained from $\CF$ by the complex conjugation   
on the coefficients. 

If a perverse sheaf $\CF$ on $P^{\nu}$ is obtained by gluing some 
irreducible factorizable sheaves into some points of $P$  
then its Verdier dual $D\CF$ is canonically isomorphic to $\bar{\CF}$. 
Therefore, the Poincar\'{e}-Verdier duality induces a perfect  
hermitian pairing on $R\Gamma(P^{\nu};\CF)$. 

Therefore,  
in notations of theorem \ref{subq},   
The Poincar\'{e}-Verdier duality 
induces a non-degenerate hermitian form on the local system 
$R^{-2n-2}\eta_*\CH^{\nu\sharp}_{\vmu}$. 

By a little more elaborated argument using fusion, one can 
introduce a canonical hermitian form on the 
systems of conformal blocks.    

Compare \cite{k}, where a certain hermitian form on the spaces 
of conformal blocks (defined up to a positive constant) has been 
introduced. 

\subsection{} 
By the similar reasons, the Verdier duality defines a hermitian form 
on all irreducible objects of $\CC$ (since the Verdier duality commutes 
with $\Phi$, cf. Theorem \ref{verdier}).

\section{Regular representation}

\subsection{} 
\label{modif} From now on we are going to modify slightly the definition 
of the categories $\CC$ and $\FS$. Let $X_{\ell}$ be the lattice 
$$
X_{\ell}=\{\mu\in X\otimes\Bbb{Q}|\ \mu\cdot Y_{\ell}\in\ell\BZ\}
$$
We have obviously $X\subset X_{\ell}$, and $X=X_{\ell}$ if $d|\ell$. 

In this part we will denote by $\CC$ a category of {\em $X_{\ell}$-graded} 
(instead of $X$-graded) finite dimensional vector spaces $M=\oplus
_{\lambda\in X_{\ell}} M_{\lambda}$ equipped with linear operators 
$\theta_i: M_{\lambda}\lra M_{\lambda-i'},\ \epsilon_i: M_{\lambda}
\lra M_{\lambda+i'}$ which satisfy the relations \ref{def c} (a), (b). 
This makes sense since $\langle d_ii,\lambda\rangle=i'\cdot\lambda\in\BZ$ 
for each $i\in I, \lambda\in X_{\ell}$. 

Also, in the definition of $\FS$ we replace $X$ by $X_{\ell}$. All the results  
of the previous parts hold true {\em verbatim} with this modification. 

We set $\dd_{\ell}=\card(X_{\ell}/Y_{\ell})$; this number is equal to 
the determinant of the form $\mu_1,\mu_2\mapsto\frac{1}{\ell}\mu_1
\cdot\mu_2$ on $Y_{\ell}$. 

\subsection{} 
\label{dot} Let $\tfu\subset\fu$ be the $k$-subalgebra generated by 
$\tK_i, \epsilon_i, \theta_i\ (i\in I)$. Following the method of 
\cite{l1} 23.1, define a new algebra $\dfu$ (without unit) as follows. 

If $\mu', \mu''\in X_{\ell}$, we set 
$$
_{\mu'}\tfu_{\mu''}=\tfu/(\sum_{i\in I} (\tK_i-\zeta^{i\cdot\mu'})\tfu+
\sum_{i\in I} \tfu(\tK_i-\zeta^{i\cdot\mu''}));\ 
\dfu=\oplus_{\mu',\mu''\in X_{\ell}} (_{\mu'}\tfu_{\mu''}).
$$
Let $\pi_{\mu',\mu''}:\tfu\lra\ _{\mu'}\tfu_{\mu''}$ be the canonical 
projection. We set $1_{\mu}=\pi_{\mu,\mu}(1)\in\dfu$. 
The structure of an algebra on $\dfu$ is defined as in {\em loc. cit.} 

As in {\em loc. cit.}, the category $\CC$ may be identified with the category 
of finite dimensional (over $k$) (left) $\dfu$-modules $M$ which are 
{\em unital}, i.e.  
 
(a) for every $x\in M$, $\sum_{\mu\in X_{\ell}}\ 1_{\mu}x=x$. 

If $M$ is 
such a module, the $X_{\ell}$-grading on $M$ is defined by $M_{\mu}=1_{\mu}M$. 

Let $\fu'$ denote  the quotient algebra of the algebra $\tfu$ by the 
relations $\tK_i^{l_i}=1\ (i\in I)$. Here $l_i:=\frac{l}{(l,d_i)}$. We have 
an isomorphism of vector spaces $\dfu=\fu'\otimes k[Y_{\ell}]$, cf. \ref{adj} 
below.   

\subsection{} Let $a:\CC\iso\CC_r$ be an equivalence defined 
by $aM=M\ (M\in\CC)$ as an $X_{\ell}$-vector space, $mx=A(x)m\ 
(x\in\fu, m\in M)$. Here $A: \fu\lra\fu^{\opp}$ is the antipode.  
We will use the same notation $a$ for a similar equivalence $\CC_r\iso\CC$. 

Let us consider the category $\CC\otimes\CC$ (resp. $\CC\otimes\CC_r$) 
which may be identified with   
the category of finite dimensional $\dfu\otimes\dfu$- (resp. 
$\dfu\otimes(\dfu)^{\opp}$-) modules satisfying 
a "unitality" condition similar to (a) above. Let us consider the 
algebra $\dfu$ itself as 
a regular $\dfu\otimes(\dfu)^{\opp}$-module. It is infinite dimensional, 
but is a union of finite 
dimensional modules, hence it may be considered as an 
object of the category $\Ind(\CC\otimes\CC_r)=\Ind\CC\otimes\Ind\CC_r$ 
where $\Ind$ denotes the category of $\Ind$-objects, cf. \cite{d4} \S 4. 
Let us denote by $\bR$ the image of this object under the equivalence 
$\Id\otimes a:\Ind\CC\otimes\Ind\CC_r\iso\Ind\CC\otimes\Ind\CC$. 

Every object $\CO\in\CC\otimes\CC$ induces a functor $F_{\CO}:\CC\lra\CC$ 
defined by 
$$
F_{\CO}(M)=a(aM\otimes_{\CC}\CO).
$$
The same formula defines a functor $F_{\CO}:\Ind\CC\lra\Ind\CC$ for 
$\CO\in\Ind(\CC\otimes\CC)$. 

We have $F_{\bR}=\Id_{\Ind\CC}$. 

We can consider a version of the above formalism using 
semiinfinite $\Tor$'s. An object $\CO\in\Ind(\CC\otimes\CC)$ defines functors 
$F_{\CO;\infh+n}:\Ind\CC\lra\Ind\CC\ (n\in\BZ)$ defined by 
$$
F_{\CO;\infh+n}(M)=a\Tor^{\CC}_{\infh+n}(aM,\CO).
$$
\subsection{Theorem} {\em {\em (i)} We have $F_{{\em\bR};\infh+n}=\Id_{\Ind\CC}$ 
if $n=0$, and $0$ otherwise. 

{\em (ii)} Conversely, suppose we have an object $Q\in\Ind(\CC\otimes\CC)$ 
together with an isomorphism of functors 
$\phi: F_{{\em\bR};\infh+\bullet}\iso F_{Q;\infh+\bullet}$. 
Then $\phi$ is induced 
by the unique isomorphism ${\em\bR}\iso Q$.} $\Box$ 

\subsection{Adjoint representation} 
\label{adj} For $\mu\in Y_{\ell}$, let $T(\mu)$ 
be a one-dimensional  
$\dfu\otimes(\dfu)^{\opp}$-module equal to $L(\mu)$ (resp. to $aL(-\mu)$) 
as a $\dfu$- (resp. $(\dfu)^{\opp}$-) module. Let us consider the module 
$T_{\mu}\bR=\bR\otimes T(\mu)\in\Ind(\CC\otimes\CC)$. This object 
represents the same functor $\Id_{\Ind\CC}$, hence we have a canonical 
isomorphism $t_{\mu}:\bR\iso T_{\mu}\bR$. 

Let us denote by $\had\in\Ind\CC$ the image of $\bR$ under the tensor 
product $\otimes:\Ind(\CC\otimes\CC)\lra\Ind\CC$. The isomorphisms $t_{\mu}$ 
above induce an action of the lattice $Y_{\ell}$ on $\had$. 
Set $\ad=\had/Y_{\ell}$. This is an object of $\CC\subset\Ind\CC$ which is 
equal to the algebra $\fu'$ considered as a $\dfu$-module by means 
of the adjoint action. 

In the notations of \ref{ms}, let us consider an object 
$$
\ad_{\MS}:=\oplus_{\mu\in\Delta_{\ell}} \langle L(\mu)\otimes L(\mu)^*
\rangle_{\MS}\in\MS,
$$
cf. \cite{bfm} 4.5.3. 
  
\subsection{Theorem} {\em We have a canonical isomorphism   
$\langle{\em\ad}\rangle_{\MS}={\em\ad}_{\MS}$.} $\Box$

\section{Regular sheaf}
\label{reg sheaf}  

\subsection{Degeneration of quadrics} The construction below is taken 
from \cite{kl}II 15.2. Let us consider the quadric $Q\subset\BP^1
\times\BP^1\times\BA^1$ given by the equation $uv=t$ where 
$(u,v,t)$ are coordinates in the triple product. Let $f:Q\lra\BA^1$ 
be the projection to the third coordinate; for $t\in\BA^1$ denote 
$Q_t:=f^{-1}(t)$. For $t\neq 0$, $Q_t$ is isomorphic to $\BP^1$; 
the fiber $Q_0$ is a union of two projective lines clutched at a point:  
$Q_0=Q_u\cup Q_v$ where $Q_u$ (resp. $Q_v$) is an irreducible component  
given (in $Q_0$) by the equation $v=0$ (resp. $u=0$) and is isomorphic 
to $\BP^1$; their intersection being a point. We set 
$'Q=f^{-1}(\BA^1-\{0\})$.  

We have two sections $x_1, x_2: \BA^1\lra Q$ given by $x_1(t)=
(\infty,0,t),\ x_2(t)=(0,\infty,t)$. Consider two "coordinate charts" 
at these points: the maps $\phi_1,\phi_2:\BP^1\times\BA^1\lra Q$ 
given by
$$
\phi_1(z,t)=(\frac{tz}{z-1},\frac{z-1}{z},t);\ 
\phi_2(z,t)=(\frac{z-1}{z},\frac{tz}{z-1},t).
$$
This defines a map 

(a) $\phi: \BA^1-\{0\}\lra\tP(2)$,  

in the notations of \ref{right}. 

\subsection{} For $\nu\in Y^+$, let us consider the corresponding 
(relative over $\BA^1$) configuration scheme  
$f^{\nu}: Q^{\nu}_{/\BA^1}\lra\BA^1$.  
For the brevity we will omit the subscript $_{/\BA^1}$ indicating that 
we are dealing with the relative version of configuration spaces. We denote 
by $Q^{\nu\bullet}$ (resp. $Q^{\nu\circ}$) the subspace of configurations 
with the points distinct from $x_1, x_2$ (resp. also pairwise distinct).  
We set $'Q^{\nu\circ}=Q^{\nu\circ}|_{\BA^1-\{0\}}$, etc. 

The map $\phi$ above, composed with the canonical projection 
$\tP(2)\lra P(2)$, induces the maps 
$$
'Q^{\nu\circ}\lra P^{\nu}(2)\ (\nu\in Y^+)
$$
(in the notations of \ref{diag}). For $\nu\in Y^+$ and 
$\mu_1,\mu_2\in X_{\ell}$ such that $\mu_1+\mu_2-\nu= 2\rho_{\ell}$, 
let $\CI_{\mu_1,\mu_2}^{\nu}$ denote the local 
system over $'Q^{\nu\circ}$ which is the inverse image of the local 
system $\CH^{\nu}_{\mu_1,\mu_2}$ over $P^{\nu}(2)^{\circ}$. 
Let $\CI_{\mu_1,\mu_2}^{\nu\bullet}$ denote the perverse sheaf 
over $'Q^{\nu\bullet}(\BC)$ which is the middle extension of 
$\CI_{\mu_1,\mu_2}^{\nu}[\dim Q^{\nu}]$. 

\subsection{}
\label{nearby}  
Let us take the nearby cycles and get a perverse sheaf $\Psi_{f^{\nu}}
(\CI^{\nu\bullet}_{\mu_1,\mu_2})$ over $Q_0^{\nu\bullet}(\BC)$. Let us consider the 
space $Q_0^{\nu\bullet}$ more attentively. This is a reducible scheme 
which is a union 
$$
Q^{\nu\bullet}_0=\bigcup_{\nu_1+\nu_2=\nu}\ \BA^{\nu_1}\times\BA^{\nu_2},
$$
the component $\BA^{\nu_1}\times\BA^{\nu_2}$ corresponding to configurations  
where $\nu_1$ (resp. $\nu_2$) points are running on the affine line 
$Q_u-x_1(0)$ (resp. $Q_v-x_2(0)$). Here we identify these affine lines 
with a "standard" one using the coordinates $u$ and $v$ respectively.    
Using this decomposition we can define a closed embedding  
$$
i_{\nu}: Q^{\nu\bullet}_0\hra\BA^{\nu}\times\BA^{\nu}
$$
whose restriction to a component $\BA^{\nu_1}\times\BA^{\nu_2}$ sends 
a configuration as above, to the configuration where all 
remaining points are equal to zero. 
Let us define a perverse sheaf 
$$
\CR_{\mu_1,\mu_2}^{\nu,\nu}=i_{\nu*}\Psi_{f^{\nu}}(\CI^{\nu\bullet}
_{\mu_2,\mu_1})
\in\CM(\BA^{\nu}(\BC)\times\BA^{\nu}(\BC);\CS)
$$
Let us consider the collection of sheaves $\{\CR_{\mu_1,\mu_2}^{\nu,\nu}|\ 
\mu_1,\mu_2\in X_{\ell}, \nu\in Y^+, \mu_1+\mu_2-\nu=2\rho_{\ell}\}$. 
One can complete this collection to an object $\CR$ of the category 
$\Ind(\FS\otimes\FS)$ where $\FS\otimes\FS$ is understood as 
a category of finite factorizable sheaves corresponding to the 
{\em square} of our initial Cartan datum, i.e. $I\coprod I$, etc. 
For a precise construction, see \cite{bfs}. 

\subsection{Theorem} {\em We have $\Phi(\CR)={\em\bR}$.} $\Box$                   
 
\newpage
\begin{center}
{\bf Part III. Modular.}
\end{center}

Almost all the results of this part are due to R.Bezrukavnikov.

\section{Heisenberg local systems} 

In this section we sketch a construction of certain 
remarkable cohesive local systems on arbitrary smooth families of  
compact smooth curves,   
to be called the {\em Heisenberg local systems}. 
 
In the definition and construction of local sustems below we will 
have to assume that our base field $k$ contains roots of unity 
of sufficiently high degree; the characteristic of $k$ is assumed to be 
prime to this degree.           

\subsection{} From now on until \ref{last heis} we fix a smooth 
proper morphism $f: C\lra S$ of relative dimension $1$,   
$S$ being a smooth connected scheme over $\BC$. For $s\in S$, we denote 
$C_s:=f^{-1}(s)$. Let $g$ be the genus of fibres of $f$.   

Let $S_{\lambda}$ denote the total space of the determinant line bundle 
$\lambda_{C/S}=\det Rf_*\Omega^1_{C/S}$ without the zero section.    
For any object (?) over $S$ (e.g., a scheme over $S$, a sheaf 
over a scheme over $S$, etc.), we will denote by (?)$_{\lambda}$ its 
base change under $S_{\lambda}\lra S$. 

Below, if we speak about a scheme as a topological (analytic)  
space, we mean its set of $\BC$-points with the usual topology 
(resp. analytic structure). 

\subsection{} 
\label{rel conf}  
We will use the relative versions of configuration spaces; to indicate 
this, we will use the subscript $_{/S}$.    
Thus, if $J$ is a finite set, $C^J_{/S}$ will denote the $J$-fold fibered 
product of $C$ with itself over $S$, etc. 

Let $C(J)_{/S}$ denote the subscheme of the $J$-fold cartesian power 
of the relative tangent bundle $T_{C/S}$ consisting of  
$J$-tuples $\{x_j, \tau_j\}$ where 
$x_j\in C$ and $\tau_j\neq 0$ in $T_{C/S,x}$, the points $x_j$ being pairwise 
distinct. 
Let $\tC(J)_{/S}$ denote the 
space of $J$-tuples of holomorphic embeddings $\phi_j: D\times S\lra C$ 
over $S$ with disjoint images; we have the $1$-jet maps 
$\tC(J)_{/S}\lra C(J)_{/S}$. 

An epimorphism $\rho: K\lra J$ induces the maps 
$$
m_{C/S}(\rho):\ \prod_J \tD(K_j)\times\tC(J)_{/S}\lra \tC(J)_{/S}
$$
and
$$
\bar{m}_{C/S}(\rho): \prod_J D(K_j)\times\tC(J)_{/S}\lra C(K)_{/S}
$$
which satisfy the compatibilities as in \ref{right}. 

\subsection{} We extend the function $n$ to $X_{\ell}$ (see \ref{modif}) by 
$n(\mu)=\frac{1}{2}\mu\cdot\mu-
\mu\cdot\rho_{\ell}\ (\mu\in X_{\ell})$. We will denote by $\CJ$ the 
$X_{\ell}$-coloured local system over the operad of disks $\CD$ which 
is defined exactly as in \ref{stand d}, with $X$ replaced by $X_{\ell}$.    

\subsection{} A {\em cohesive local system $\CH$ of level $\mu\in X_{\ell}$ 
over $C/S$} 
is a collection of local systems $\CH(\pi)$ over the spaces 
$C(J)_{/S;\lambda}$ 
given for every map $\pi: J\lra X_{\ell}$ of level $\mu$ (note the base change  
to $S_{\lambda}$!), together with the factorization isomorphisms
$$
\phi_C(\rho): m_{C/S}(\rho)^*\CH(\pi)\iso\Boxtimes_J \CJ(\pi_j)
\boxtimes \CH(\rho_*\pi).
$$
Here we have denoted by the same letter $\CH(\pi)$ the lifting of $\CH(\pi)$ 
to $\tC_{/S;\lambda}$. 
The factorization isomorphisms must satisfy the obvious analogs 
of properties \ref{cls def} (a), (b). 

\subsection{} 
\label{start heis} Now we will sketch a construction of certain cohesive 
local system over $C/S$ of level $(2-2g)\rho_{\ell}$. For alternative 
beautiful constructions of $\CH$, see \cite{bp}.   

To simplify the exposition we will assume below that $g\geq 2$ 
(the construction for $g\leq 1$ needs some modification, and we omit it here, 
see \cite{bfs}).   
Let us consider the group scheme $\Pic(C/S)\otimes X_{\ell}$ over $S$. 
Here $\Pic(C/S)$ is the relative Picard scheme. The group 
of connected components $\pi_0(\Pic(C/S)\otimes X_{\ell})$ is equal 
to $X_{\ell}$. Let us denote by $\Jac$  
the connected component corresponding to the element $(2-2g)\rho_{\ell}$; 
this is an 
abelian scheme over $S$, due to the existence of the section 
$S\lra\Jac$ defined by $\Omega^1_{C/S}\otimes(-\rho_{\ell})$.  

For a scheme $S'$ over $S$, let $H_1(S'/S)$ denote the local system 
of the first relative integral homology groups over $S$. We have 
$H_1(\Jac/S)=H_1(C/S)\otimes X_{\ell}$. We will denote by $\omega$ the 
polarization of $\Jac$ (i.e. the  
skew symmetric form on the latter local system) equal to the tensor  
product of the standard form on $H_1(C/S)$ and the form 
$(\mu_1,\mu_2)\mapsto\frac{\dd^g_{\ell}}{\ell}\mu_1\cdot\mu_2$ on $X_{\ell}$. 
Note that the assumption $g\geq 2$ implies that 
$\frac{\dd^g_{\ell}}{\ell}\mu_1\cdot\mu_2\in\BZ$ for any $\mu_1,\mu_2\in 
X_{\ell}$.   
Since the latter form is positive definite, $\omega$ is relatively ample 
(i.e. defines a relatively ample invertible sheaf on $\Jac$). 

\subsection{} Let $\alpha=\sum n_{\mu}\cdot\mu\in\BN[X_{\ell}]$; set 
$\Supp(\alpha)=\{\mu|\ n_{\mu}\neq 0\}$. Let us say that $\alpha$ 
is {\em admissible} if $\sum n_{\mu}\mu=(2-2g)\rho_{\ell}$. 
Let us denote by 
$$
\ya_{\alpha}: C^{\alpha}_{/S}\lra\Pic(C/S)\otimes X_{\ell}
$$
the Abel-Jacobi map sending $\sum \mu\cdot x_{\mu}$ to 
$\sum x_{\mu}\otimes\mu$. If $\alpha$ is admissible then the map 
$\ya_{\alpha}$ lands in $\Jac$. 

Let $D^{\alpha}$ denote the following relative divisor on $C^{\alpha}_{/S}$  
$$
D^{\alpha}=\frac{\dd_{\ell}^g}{\ell}(\sum_{\mu\neq\nu}\mu\cdot\nu
\Delta_{\mu\nu}+\frac{1}{2}\sum_{\mu}\mu\cdot\mu\Delta_{\mu\mu}).
$$
Here $\Delta_{\mu\nu}\ (\mu,\nu\in \Supp(\alpha))$ denotes  
the corresponding diagonal in $C^{\alpha}_{/S}$. Note that all the 
multiplicities are integers.  

Let $\pi: J\lra X_{\ell}$ be an unfolding of $\alpha$. We will denote by 
$D^{\pi}$ the pull-back of $D^{\alpha}$ to $C^J_{/S}$. Let us 
introduce the following line bundles  
$$
\CL(\pi)=\otimes_{j\in J}\CT_j^{\otimes\frac{\dd^g_{\ell}}{\ell}n(\pi(j))}
\otimes\CO(D^{\pi})
$$
on $C^J_{/S}$, and 
$$
\CL_{\alpha}=\CL(\pi)/\Sigma_{\pi}
$$
on $C^{\alpha}_{/S}$ (the action of $\Sigma_{\pi}$ is an obvious one). 
Here $\CT_j$ denotes the relative tangent line bundle on $C^J_{/S}$ in the 
direction $j$. Note that the numbers $\frac{\dd^g_{\ell}}{\ell}n(\mu)\ 
(\mu\in X_{\ell})$ are integers.  

\subsection{Proposition} {\em There exists a unique line bundle 
$\CL$ on $\Jac$ such that for each admissible $\alpha$, we have 
$\CL_{\alpha}=\ya_{\alpha}^*(\CL)$. The first Chern class 
$c_1(\CL)=-[\omega]$.} $\Box$ 

\subsection{} In the sequel if $\CL_0$ is a line bundle, let 
$\dCL_0$ denote the total space its with the zero section removed. 
 
The next step is the construction of a certain local system $\fH$ over 
$\dCL_{\lambda}$. Its dimension is equal to $\dd^g_{\ell}$ and 
the monodromy 
around the zero section of $\CL$ (resp. of the determinant bundle) is equal to 
$\zeta^{-2\ell/\dd_{\ell}^g}$ (resp. $(-1)^{\rk(X)}\zeta^{-12\rho_{\ell}\cdot
\rho_{\ell}})$. The construction of $\fH$ is outlined below. 

The previous construction assigns to a triple 

(a lattice $\Lambda$, a symmetric bilinear form $(\ ,\ ):\ \Lambda\times\Lambda
\lra\BZ$, $\nu\in\Lambda$) 

an abelian scheme $\Jac_{\Lambda}:=(\Pic(C/S)\otimes\Lambda)_{(2g-2)\nu}$ 
over $S$, together with a line bundle $\CL_{\Lambda}$ on it (in the 
definition of $\CL_{\Lambda}$ one should use the function 
$n_{\nu}(\mu)=\frac{1}{2}\mu\cdot\mu+\mu\cdot\nu$). We considered the case 
$\Lambda=X_{\ell},\ (\mu_1,\mu_2)=\frac{\dd^g_{\ell}}{\ell}\mu_1\cdot\mu_2,\ 
\nu=-\rho_{\ell}$.  

Now let us apply this construction to the lattice 
$\Lambda=X_{\ell}\oplus Y_{\ell}$, the bilinear form 
$((\mu_1,\nu_1),(\mu_2,\nu_2))=-\frac{1}{\ell}(\nu_1\cdot\mu_2+\nu_2\cdot\mu_1
+\nu_1\cdot\nu_2)$ and $\nu=(-\rho_{\ell},0)$. The first projection 
$\Lambda\lra X_{\ell}$ induces the morphism 
$$
p:\ \Jac_{\Lambda}\lra\Jac
$$
the fibers of $p$ are abelian varieties $\Jac(C_s)\otimes Y_{\ell}\ (s\in S)$.  

\subsection{Theorem} {\em {\em (i)} The line bundle $\CL_{\Lambda}$ 
is relatively ample with respect to $p$. The direct image 
$\CE:=p_*\CL_{\Lambda}$ is  a locally free sheaf of rank $\dd_{\ell}^g$. 

{\em (ii)} We have an isomorphism
$$
\det(\CE)=\CL\otimes\lambda^{\dd^g_{\ell}(-\frac{1}{2}\rk(X_{\ell})+
6\frac{\rho_{\ell}\cdot\rho_{\ell}}{\ell})}.
$$} $\Box$ 

Here $\lambda$ denotes the pull-back of the determinant bundle 
$\lambda_{C/S}$ to $\Jac$. 

\subsection{} Let us assume for a moment that $k=\BC$ and 
$\zeta=\exp(-\frac{\pi\sqrt{-1}}{\ell})$.   
By the result of Beilinson-Kazhdan, \cite{bk} 4.2, 
the vector bundle $\CE$ carries a canonical flat projective connection. 
By {\em loc. cit.} 2.5, its lifting to $\det(\CE)^{\bullet}$ carries a flat 
connection with the scalar monodromy around the zero section equal to 
$\exp(\frac{2\pi\sqrt{-1}}{\dd^g_{\ell}})$. We have an 
obvious map 
$$
m: \dCL_{\lambda}\lra\CL\otimes\lambda^
{\dd^g_{\ell}(-\frac{1}{2}\rk(X_{\ell})+
6\frac{\rho_{\ell}\cdot\rho_{\ell}}{\ell})}.
$$
By definition, $\fH$ is the local system of horizontal sections 
of the pull-back of $\CE$ to $\dCL_{\lambda}$. The claim about its monodromies  
follows from part (ii) of the previous theorem. 

This completes the construction of $\fH$ for $k=\BC$ and $\zeta=\exp
(-\frac{\pi\sqrt{-1}}{\ell})$. The case of arbitrary $k$ (of sufficiently 
large characteristic) and $\zeta$ follows from this one.       

\subsection{} 
\label{last heis}  
Let us consider an obvious map $q: C(J)_{/S;\lambda}\lra C^J_{/S;\lambda}$.  
The pull-back  
$q^*\CL(\pi)$ has a canonical non-zero section $s$. Let $\tCH(\pi)$ be 
the pull-back of the local system $\fH$ to $q^*\CL(\pi)$. By definition, 
we set $\CH(\pi)=s^*\tCH(\pi)$. For the construction of the factorization 
isomorphisms, see \cite{bfs}.  

\subsection{} 
\label{level} Let $\fg$ be the simple Lie algebra connected with 
our Cartan datum.   
Assume that $\zeta=\exp(\frac{\pi\sqrt{-1}}{d\kappa})$ for some 
positive integer $\kappa$, cf. \ref{charge} ($d$ is defined in \ref{root}).   

We have $12\rho_{\ell}\cdot\rho_{\ell}\equiv 12\rho\cdot\rho\ \modul\ l$, and 
$\rk(X)\equiv\dim\fg\ \modul\ 2$. By the strange formula of 
Freudenthal-de Vries, we have $12\rho\cdot\rho=dh\dim\fg$ where $h$ 
is the dual Coxeter number of our Cartan datum.   
It follows that the monodromy of 
$\CH$ around the zero section of the determinant line bundle is equal to 
$\exp(\pi\sqrt{-1}\frac{(\kappa-h)\dim\fg}{\kappa})$. This number coincides 
with the multiplicative central charge of the conformal 
field theory associated with the affine Lie algebra $\hfg$ at level $\kappa$ 
(see \cite{bfm} 4.4.1, 6.1.1, 2.1.3, \cite{tuy} 1.2.2), cf 
\ref{subquot} below. 

{\em UNIVERSAL HEISENBERG SYSTEMS} 

\subsection{} 
Let us define a category $\Sew$ as follows (cf. \cite{bfm} 4.3.2). 
Its object $A$ is a finite set $\bA$ together with a collection $N_A=
\{n\}$ of non-intersecting two-element subsets $n\subset\bA$. Given 
such an object, we set $A^1=\bigcup_{N_A}\ n,\ A^0=\bA-A^1$. A morphism 
$f: A\lra B$ is an embedding $i_f: \baB\hra \bA$ and a collection $N_f$ of 
non-intersecting two-element subsets of $\baB-\bA$ such that 
$N_A=N_B\coprod N_f$. The composition of morphisms is obvious. ($\Sew$  
coincides with the category $Sets^{\sharp}/\emp$, in the notations of 
\cite{bfm} 4.3.2.)  

For $A\in\Sew$, let us call an {\em $A$-curve} a data $(C, \{x_a,\tau_a\}
_{A^0})$ where $C$ is a smooth proper (possibly disconnected) complex  
curve, $\{x_a,\tau_a\}_{A^0}$ is an $A^0$-tuple of distinct points 
$x_a\in C$ together with non-zero tangent vectors $\tau_a$ at them.  
For such a curve, let $\bC_A$ denote the curve obtained from $C$ by clutching  
pairwise the points $x_{a'}$ with $x_{a''}\ (n=\{a',a''\})$ 
for all sets $n\in N_A$. Thus, the set $N_A$ is in the bijection with 
the set of nodes of the curve $\bC_A$ ($\bC_A=C$ if $N_A=\emp$). 

Let us call an {\em enhanced graph} a pair $\Gamma=(\bGamma,\bg)$ here $\bGamma$ 
is a non-oriented graph and $\bg=\{g_v\}_{v\in\Ve(\bGamma)}$ is a $\BN$-
valued $0$-chain of $\bGamma$. Here $\Ve(\bGamma)$ denotes the set of 
vertices of $\bGamma$. Let us assign to a curve $\bC_A$ an enhanced graph 
$\Gamma(\bC_A)=(\bGamma(\bC_A), \bg(\bC_A))$. By definition,  
$\bGamma(\bC_A)$ is a graph with $\Ve(\bGamma(\bC_A))=\pi_0(C)=\{$ 
the set of irreducible components of $\bC_A\}$ and the set of edges   
$\Ed(\bGamma(\bC_A))=N_A$, an edge $n=\{a',a''\}$ connecting the 
vertices corresponding 
to the components of the points $x_{a'}, x_{a''}$. For $v\in\pi_0(C)$,  
$g(\bC_A)_v$ is equal to the genus of the corresponding component 
$C_v\subset C$.   

\subsection{}   
Let $\CM_A$ denote the moduli stack of $A$-curves $(C,\ldots)$ such that 
the curve $\bC_A$ is stable in the sense of \cite{dm} (in particular 
connected). The stack $\CM_A$ is smooth; we have   
$\CM_A=\coprod_{g\geq 0}\CM_{A,g}$ where $\CM_{A,g}$ is a  
substack of $A$-curves $C$ with $\bC_A$ having genus (i.e. 
$\dim H^1(\bC_A,\CO_{\bC_A})$) equal to $g$. In turn, we have the decomposition 
into connected components 
$$
\CM_{A,g}=\coprod_{\Gamma,\vA^0}\ \CM_{\vA^0,g,\Gamma}
$$ 
where $\CM_{\vA^0,g,\Gamma}$ is the stack of $A$-curves $(C,\{x_a,\tau_a\})$ as 
above, with $\Gamma(\bC_A)=\Gamma$, $\vA^0=\{A^0_v\}_{v\in\Ve(\bGamma)}$, 
$A^0=\coprod\ A^0_v$, such that $x_a$ lives on the connected component 
$C_v$ for $a\in A^0_v$.   
  
We denote by $\eta: C_{A,g}\lra\CM_{A,g}$ 
(resp. $\baeta: \bC_{A,g}\lra\CM_{A,g})$ the universal smooth curve 
(resp. stable surve).   
For $\nu\in Y^+$, we have the corresponding {\em relative} configuration 
spaces $C^{\nu}_{A,g}, C^{\nu\circ}_{A,g}, \bC^{\nu}_{A,g}$. For brevity we 
omit the relativeness subscript $/\CM_{A,g}$ from these notations. 
The notation $C^{\nu}_{\vA^0,g,\Gamma}$ etc., will mean the restriction 
of these 
configuration spaces to the component $\CM_{\vA^0,g,\Gamma}$. 

Let $\CM_g$ be the moduli stack of smooth connected curves of genus $g$, 
and $\bCM_g$ be its Grothendieck-Deligne-Mumford-Knudsen  
compactification, i.e. the moduli stack of stable curves of genus $g$. 
Let $\baeta: \bC_g\lra\bCM_g$ be the universal stable curve; let 
$\blambda_g=\det R\baeta_*(\omega_{\bC_g/\bCM_g})$ be the determinant line 
bundle; let $\bCM_{g;\lambda}\lra\bCM_g$ be its total space with the zero 
section removed.      

We have obvious maps $\CM_{A,g}\lra\bCM_g$. 
Let the complementary subscript $(\cdot)_{\lambda}$ denote the base change of 
all the above objects under $\bCM_{g;\lambda}\lra\bCM_g$. 

\subsection{} 
\label{def h} Let us consider the configuration space $C^{\nu\circ}_{A,g}$; 
it is the moduli stack of $\nu$ distinct points running on $A$-curves 
$(C,\{x_a,\tau_a\})$ and not equal to the marked points $x_a$. This stack 
decomposes into connected components as follows: 
$$
C^{\nu\circ}_{A,g}=\coprod_{\Gamma,\vA^0,\vnu}\ C^{\vnu}_{\vA^0,g,\Gamma}
$$
where $\vnu=\{\nu_v\}_{v\in\Ve(\Gamma)}$ and $C^{\vnu}_{\vA^0,g,\Gamma}$ being 
the moduli stack of objects as above, with $\Gamma(\bC_A)=\Gamma$ and 
$\nu_v$ points running on the component $C_v$. The decomposition is taken 
over appropriate graphs $\Gamma$, decompositions $A^0=\coprod A^0_v$ and 
the tuples $\vnu$ with $\sum \nu_v=\nu$. 

Let us call an $A_0$-tuple $\vmu=\{\mu_a\}\in X_{\ell}^{A^0}$ 
{\em $(g,\nu)$-good} if 

(a) $\sum_{a\in A^0}\ \mu_a-\nu\equiv (2-2g)\rho_{\ell}\ \modul\ Y_{\ell}$.

Given such a tuple, we are going to define certain local system 
$\CH^{\nu}_{\vmu;A,g}$ over $C^{\nu\circ}_{A,g;\lambda}$. Let us describe its 
restriction $\CH^{\vnu}_{\vmu;\vA^0,g,\Gamma}$ to a connected 
component $C^{\vnu}_{\vA^0,g,\Gamma;\lambda}$. 

Let $\Gamma'$ be the first subdivision of $\bGamma$. We have 
$\Ve(\Gamma')=\Ve(\bGamma)\coprod\Ed(\bGamma)=\pi_0(C)\coprod N_A$. 
The edges of $\Gamma'$ are 
indexed by the pairs $(n,a)$ where $n\in N_A, a\in n$, the corresponding 
edge $e_{n,a}$ having the ends $a$ and $n$. Let us define an orientation 
of $\Gamma'$ by the requierement that $a$ is the beginning of $e_{n,a}$. 
Consider the chain complex 
$$
C_1(\Gamma'; X_{\ell}/Y_{\ell})\overset{d}{\lra}C_0(\Gamma';X_{\ell}/Y_{\ell}). 
$$
Let us define a $0$-chain $c=c^{\vnu}_{\vmu}\in C_0(\Gamma';X_{\ell}/Y_{\ell})$ 
by 
$$
c(v)=\sum_{a\in A^0_v}\ \mu_a+(2g_v-2)\rho_{\ell}-\nu_v\ (v\in\pi_0(C));\ 
c(n)=2\rho_{\ell}\ (n\in N_A). 
$$
The goodness assumption (a) ensures that $c$ is a boundary.      
By definition, 
$$
\CH^{\vnu}_{\vmu;\vA^0,g,\Gamma}=\oplus_{\chi:\ d\chi=c}\ \CH_{\chi}
$$
Note that the set $\{\chi|\ d\chi=c\}$ is a torsor over the group 
$H_1(\Gamma';X_{\ell}/Y_{\ell})=H_1(\Gamma;X_{\ell}/Y_{\ell})$.
 
The local system $\CH_{\chi}$ is defined below, in \ref{h chi}, after a little 
notational preparation. 

\subsection{} Given two finite 
sets $J, K$, let $C(J;K)_g$ denote the moduli stack of objects 
$(C,\{x_j\},\{y_k,\tau_k\})$. Here $C$ is a smooth proper connected 
curve of genus $g$, $\{x_j\}$ is a $J$-tuple of distinct points $x_j\in C$ 
and $\{y_k,\tau_k\}$ is a $K$-tuple of distinct points $y_k\in C$ together 
with non-zero tangent vectors $\tau_k\in T_{y_k}C$. We suppose that 
$y_k\neq x_j$ for all $k,j$. 

We set $C(J)_g:=C(\emp;J)_g$. We have the forgetful maps 
$C(J\coprod K)_g\lra C(J;K)_g$. 

The construction of \ref{start heis} --- \ref{last heis} defines 
the Heisenberg system $\CH(\pi)$ over the smooth stack  
$C(J)_{g;\lambda}$ for each $\pi: J\lra X_{\ell}$.   

Given $\nu\in Y^+$, choose an unfolding of $\nu$, $\pi: J\lra I$, and set 
$C^{\nu}(K)^{\circ}_g:=C(J;K)/\Sigma_{\pi}$. Given a $K$-tuple $\vmu=
\{\mu_k\}\in X_{\ell}^K$, define a map $\tpi: J\coprod K\lra X_{\ell}$ by 
$\tpi(j)=-\pi(j)\in - I\subset X_{\ell}\ (j\in J),\ 
\tpi(k)=\mu_k\ (k\in K)$. The local system $\CH(\tpi)$ over 
$C(J\coprod K)_{g;\lambda}$ descends to $C(J;K)_{g;\lambda}$ since 
$\zeta^{2n(-i)}=1$, and then to $C^{\nu}(K)_{g;\lambda}^{\circ}$, by 
$\Sigma_{\pi}$-equivariance; denote the latter local system by 
$\tCH_{\vmu}^{\nu}$, and set $\CH^{\nu}_{\vmu}=\tCH^{\nu}_{\vmu}\otimes 
\Sign^{\nu}$, cf. \ref{sign p}. 

\subsection{Lemma} {\em If  
$\vmu\equiv\vmu'\ \modul\ Y_{\ell}^K$ then we   
have canonical isomorphisms $\CH^{\nu}_{\vmu}=\CH^{\nu}_{\vmu'}$.} $\Box$ 

Therefore, it makes sense to speak about $\CH^{\nu}_{\vmu}$ for 
$\vmu\in (X_{\ell}/Y_{\ell})^K$. 

\subsection{} 
\label{h chi} Let us return to the situation at the end of \ref{def h}. 
We have $\Gamma=(\bGamma,\{g_v\}_{v\in\Ve(\bGamma)})$. Recall that 
$A^0=\coprod\ A^0_v$; we have also $A^1=\coprod A^1_v$ where $A^1_v:= 
\{a\in A^1|\ x_a\in C_v\}\ (v\in \Ve(\bGamma)$. Set $\bA_v=A^0_v\coprod A^1_v$, 
so that $\bA=\coprod \bA_v$. We have an obvious map 

(a) $\CM_{\vA^0,g,\Gamma}\lra\prod_{\Ve(\bGamma)}\ C(\bA_v)_{g_v}$, 

and a map 

(b) $C^{\vnu}_{\vA^0,g,\Gamma}\lra \prod_{\Ve(\bGamma)}\ 
C^{\nu_v}(\bA_v)^{\circ}_{g_v}$ 

over (a). For each $v$, define an $\bA_v$-tuple $\vmu(\chi;v)$ equal to $\mu_a$ 
at $a\in A^0_v$ and $\chi(e_{a,n})$ at $a\in n\subset A^1$. By definition, 
the local system $\CH_{\chi}$ over $C^{\vnu}_{\vA^0,g,\Gamma;\lambda}$ is the 
inverse image of the product 

$\Boxtimes_{\Ve(\bGamma)}\ \CH^{\nu_v}_{\vmu(\chi;v)}$ 

under the map (b) (pulled back to the determinant bundle).  

This completes the definition of the local systems $\CH^{\nu}_{\vmu;A,g}$. 
They have a remarkable compatibility property (when the object $A$ varies)  
which we are going to describe below, see theorem \ref{cart}.  

\subsection{} 
\label{conf fib} Let $\tCT^{\nu}_{A,g}$ denote the fundamental groupoid 
$\pi(C^{\nu\circ}_{A,g;\lambda})$. We are going to show that these groupoids 
form a cofibered category over $\Sew$. 

\subsection{} 
\label{cofib} A morphism $f: A\lra B$ in $\Sew$ is called  
a {\em sewing} (resp. {\em deleting}) if $A^0=B^0$ (resp.  
$N_f=\emp$). A sewing $f$ with $\card(N_f)=1$ is called {\em simple}.  
Each morphism is a composition of a sewing and a deleting; each sewing 
is a composition of simple ones.  

(a) Let $f: A\lra B$ be a simple sewing. We have canonical morphisms 
$$
\wp_f: \CM_{A,g;\lambda}\lra\dT_{\dpar\bCM_{B,g;\lambda}}\bCM_{B,g;\lambda}
$$
and 
$$
\wp_f^{(\nu\circ)}:C^{\nu\circ}_{A,g;\lambda}\lra
\dT_{\dpar\overline{C^{\nu\circ}_{B,g;\lambda}}}\ \overline
{C^{\nu\circ}_{B,g;\lambda}},
$$
over $\wp_f$, cf. \cite{bfm} 4.3.1. Here $\bCM_{B,g}$ 
(resp. $\overline{C^{\nu\circ}_{B,g}}$)  
denotes the Grothendieck-Deligne-Mumford-Knudsen 
compactification of $\CM_{B,g}$ (resp. of $C^{\nu\circ}_{B,g}$)  
and $\dpar\bCM_{B,g}$ (resp. $\dpar\overline{C^{\nu\circ}_{B,g}}$)  
denotes the smooth locus of the boundary $\bCM_{B,g}-\CM_{B,g}$ (resp. 
$\overline{C^{\nu\circ}_{B,g}}-\bC^{\nu}_{B,g}$).    
The subscript $_{\lambda}$ indicates the base change to the 
determinant bundle, as before. 

Composing the specialization with the inverse image under $\wp_f^{(\nu\circ)}$, 
we get the canonical map $f_*:\tT^{\nu}_{A,g}\lra\tT^{\nu}_{B,g}$. 

(b) Let $f:A\lra B$ be a deleting. It induces 
the obvious morphisms (denoted by the same letter)
$$
f:\CM_{A,g;\lambda}\lra\CM_{B,g;\lambda}
$$
and
$$
f:C^{\nu\circ}_{A,g;\lambda}\lra C^{\nu\circ}_{B,g;\lambda}.
$$
The last map induces $f_*:\tT^{\nu}_{A,g}\lra\tT^{\nu}_{B,g}$. 

Combining the  constructions (a) and (b) above, we get a category 
$\tT^{\nu}_g$  
cofibered in groupoids over $\Sew$, with fibers $\tT^{\nu}_{A,g}$. 

\subsection{}
\label{rep}  
Let $\Rep^{\nu}_{c;A,g}$ be the category of finite dimensional representations 
of $\tT^{\nu}_{A,g}$ (over $k$) with the monodromy $c\in k^*$ around 
the zero section of the determinant bundle. The previous construction 
shows that these categories form a fibered category $\Rep^{\nu}_{c;g}$ over 
$\Sew$. 

\subsection{} For $A\in\Sew$, let us call an $A^0$-tuple 
$\vmu=\{\mu_a\}\in X^{A^0}_{\ell}$  
{\em good} if $\sum_{A^0}\ \mu_a\in Y$. 

If $f: B\lra A$ is a morphism, define a $B^0$-tuple $f^*\vmu=\{\mu'_b\}$ by 
$\mu'_b=\mu_{i_f^{-1}(b)}$ if $b\in i_f(A^0)$, and $0$ otherwise. 
Obviously, $\sum_{B^0}\ \mu'_b=\sum_{A^0}\ \mu_a$. 

Given a good $\vmu$, let us pick an element $\nu\in Y^+$ such that 
$\nu\equiv\sum\mu_a+(2g-2)\rho_{\ell}\ \modul\ Y_{\ell}$. We    
can consider the local system $\CH^{\nu}_{\vmu;A,g}$   
as an object of $\Rep^{\nu}_{c;A,g}$ where 
$c=(-1)^{\card(I)}\zeta^{-12\rho\cdot\rho}$. 

\subsection{Theorem} 
\label{cart} {\em For any morphism $f: B\lra A$ in $\Sew$ and 
a $g$-good $\vmu\in X^{A^0}_{\ell}$, we have the canonical isomorphism 
$$
f^*\CH_{\vmu;A,g}^{\nu}=\CH_{f^*\vmu;B,g}^{\nu}.
$$
In other words, the local systems $\CH_{\vmu;B,g}^{\nu}$ define a 
{\em{\bf cartesian section}} 
of the fibered category $\Rep^{\nu}_{c;g}$ over $\Sew/A$. Here 
$c=(-1)^{\card(I)}\zeta^{-12\rho\cdot\rho}$. $\Box$}

\section{Fusion structures on $\FS$}
\label{fusion} 

\subsection{} Below we will construct a family of "fusion structures"  
on the category $\FS$ (and hence, due to the equivalence 
$\Phi$, on the category $\CC$) indexed by $m\in\BZ$. We should explain 
what a fusion structure is. This is done in \ref{def fus} below.  
We will use a modification of the formalism from \cite{bfm} 4.5.4.  

\subsection{}  
Recall that we have a regular object $\CR\in\Ind
(\FS^{\otimes 2})$, cf. Section \ref{reg sheaf}. We have the canonical 
isomorphism $t(\CR)=\CR$ 
where $t: \Ind(\FS^{\otimes 2})\iso\Ind(\FS^{\otimes 2})$ is the permutation, 
hence an object $\CR_n\in\Ind(\FS^{\otimes 2})$ is well defined for 
any two-element set $n$. 

For an object $A\in\Sew$, we set $\tA=A^0\coprod N_A$. 
Let us call an {\em $A$-collection} of factorizable sheaves an 
$\tA$-tuple $\{\CX_{\ta}\}_{\ta\in\tA}$ where $\CX_{\ta}\in\FS_{c_{\ta}}$ if 
$\ta\in A^0$ and 
$\CX_{\ta}=\CR_{\ta}$ if $\ta\in N_A$. We impose the condition that 
$\sum_{a\in A^0}\ c_a=0\in X_{\ell}/Y$.    
We will denote such an object  
$\{\CX_a;\CR_n\}_A$. It defines an object 
$$
\otimes_A\ \{\CX_a; \CR_n\}\ :=(\otimes_{a\in A^0}\ \CX_a)\otimes(\otimes
_{n\in N_A}\ \CR_n)\in\Ind(\FS^{\otimes\bA}).
$$
If $f: B\lra A$ is a morphism in $\Sew$, we define a $B$-collection 
$f^*\{\CX_a;\CR_n\}_B=\{\CY_{\tb}\}_{\tB}$ by $\CY_{\tb}=\CX_{i_f^{-1}(\tb)}$ 
for $\tb\in i_f(A^0)$, $\One$ if $\tb\in B^0-i_f(A^0)$ and $\CR_{\tb}$ if 
$\tb\in N_B=N_A\cup N_f$.    

\subsection{}    
Given an $A$-collection $\{\CX_a,\CR_n\}_A$ with $\lambda(\CX_a)=\mu_a$ 
such that 

(a) $\nu:=\sum_{A^0}\ \mu_a+(2g-2)\rho_{\ell}\in Y^+$,  

one constructs (following the pattern of \ref{glu thm}) a perverse sheaf  
$\Boxtimes_{A,g}^{(\nu)}\ \{\CX_a; \CR_n\}$ over 
$\bC^{\nu}_{A,g;\lambda}$. 
It is obtained by planting factorizable sheaves $\CX_a$ into the 
universal sections $x_a$ of the stable curve $\bC_{A,g}$, the regular 
sheaves $\CR_{n}$ into the nodes $n$ of this curve and pasting them together  
into one sheaf by the Heisenberg system $\CH_{\vmu;A,g}^{\nu}$.  

\subsection{}  
Given $A\in\Sew$ and an $A$-collection $\{\CX_a;\CR_n\}_A$,   
choose elements $\mu_a\geq\lambda(\CX_a)$ in $X_{\ell}$ such that
(a) above holds (note that $2\rho_{\ell}\in Y$). Below $\nu$ will denote the 
element as in (a) above.  

Let $\CX'_a$ denote the factorizable sheaf 
isomorphic to $\CX_a$ obtained from it by the change of the highest 
weight from $\lambda(\CX_a)$ to $\mu_a$.  
For each $m\in\BZ$, define a local system 
$\langle\otimes_A\ \{\CX_a; \CR_n\}\rangle^{(m)}_g$  
over $\CM_{A,g;\lambda}$ as follows. 

Let $\langle\otimes_A\ \{\CX_a; \CR_n\}\rangle^{(m)}
_{\vA^0,g,\Gamma}$ denote the restriction of the local system to be defined to 
the connected component $\CM_{\vA^0,g,\Gamma;\lambda}$. By definition,  
$$
\langle\otimes_{A}\ \{\CX_a; \CR_n\}\rangle^{(m)}_{\vA^0,g,\Gamma}:=
R^{m-\dim \CM_{\vA^0,g,\Gamma;\lambda}}
\baeta_*^{\nu}(\Boxtimes_{\vA^0,g,\Gamma}^{(\nu)}\ \{\CX'_a;\CR_n\}).
$$
Here 
$\Boxtimes_{\vA^0,g,\Gamma}^{(\nu)}\ \{\CX'_a;\CR_n\}$ denotes the perverse 
sheaf  
$\Boxtimes_{A,g}^{(\nu)}\ \{\CX'_a;\CR_n\}$ restricted to the subspace 
$\bC^{\nu}_{\vA^0,g,\Gamma;\lambda}$.   
This definition does not depend on the choice of the elements $\mu_a$.  

\subsection{} Given a morphism $f: B\lra A$, we define, acting as in 
\ref{conf fib}, \ref{rep}, a perverse sheaf 
$f^*(\Boxtimes^{(\nu)}_{A,g}\ \{\CX_a;\CR_n\})$ over 
$\bC^{\nu}_{B,g;\lambda}$ and local systems 
$f^*\langle\otimes_A\ \{\CX_a;\CR_n\}\rangle_g^{(m)}$ over 
$\CM_{B,g;\lambda}$.  

\subsection{Theorem} {\em In the above notations, we have 
canonical isomorphisms 
$$
f^*(\Boxtimes^{(\nu)}_{A,g}\ \{\CX_a; \CR_n\})=\Boxtimes_{B,g}^{(\nu)}\ 
f^*\{\CX_a; \CR_n\}.\ 
$$}

This is a consequence of Theorem \ref{cart} above and the definition of 
the regular sheaf $\CR$ as a sheaf of nearby cycles of the braiding local 
system, \ref{nearby}. 

\subsection{Corollary} 
\label{fus cor} {\em We have canonical isomorphisms of local 
systems  
$$
f^*\langle\otimes_A\ \{\CX_a; \CR_n\}\rangle^{(m)}_g=\langle\otimes_B\ 
f^*\{\CX_a; \CR_n\}
\rangle^{(m)}_g\ (m\in\BZ).\ \Box
$$}
\subsection{} 
\label{def fus} The previous corollary may be expressed as follows. The 
various $A$-collections of 
factorizable sheaves (resp. categories $\tRep_{c;A}$ of finite dimensional 
representations of "Teichm\"{u}ller groupoids" 
$\tTeich_{A}=\pi(\CM_{A;\lambda})$ having monodromy $c$ around the 
zero section of the determinant bundle)  
define a fibered category $\FS^{\sharp}$ (resp. $\tRep_{c}$) over $\Sew$.

For any $m\in \BZ$, the collection of local systems 
$\langle\otimes_A\ \{\CX_a; \CR_n\}\rangle^{(m)}_g$ and 
the canonical isomorphisms of the previous theorem define a 
{\bf cartesian functor}
$$
\langle\ \rangle^{(m)}:\FS^{\sharp}\lra\tRep_{c}  
$$
where $c=(-1)^{\card(I)}\zeta^{-12\rho\cdot\rho}$. We call such a functor 
a {\em fusion structure of multiplicative central 
charge $c$} on the category $\FS$.   
The category $\FS$ with this fusion structure 
will be denoted by $\FS^{(m)}$. 

The difference from the definition of a fusion structure given in 
\cite{bfm} 4.5.4 is that  
our fibered categories live over $\Sew=\Sets^{\sharp}/\emp$ and not 
over $\Sets^{\sharp}$, as in {\em op. cit.}   
 
\subsection{Example} For $A\in\Sew$, let us consider an $A$-curve 
$P=(\BP^1,\{x_a,\tau_a\})$; it defines a geometric point $Q$  
of the stack $\CM_{A,g}$ where $g=\card(N_A)$ and hence a geometric point 
$P=(Q,1)$ of the stack $\CM_{A,g;\lambda}$ since the determinant bundle is 
canonically trivialized at $Q$. 

\subsubsection{} {\bf Theorem.} {\em For an $A$-collection 
$\{\CX_a; \CR_n\}$, the stalk of the local system 
$\langle\otimes_A\ \CX_a\rangle^{(m)}_g$ at a point $P$ is isomorphic to 
$$
\Tor_{\infh-m}^{\CC}({\em \One}_r,(\otimes_{A^0}\ \Phi(\CX_a))\otimes
{\em \ad}^{\otimes g}).
$$}
To prove this, one should apply theorem \ref{integral} and the following 
remark. Degenerating all nodes into cusps, one can include the nodal 
curve $\overline{P}_A$ into a  
one-parameter family whose special fiber is $\BP^1$, with 
$\card(A^0)+g$ marked points. The nearby cycles of the sheaf 
$\Boxtimes_A\ \{\CX_a; \CR_n\}$ will be the sheaf obtained by the gluing of the 
sheaves $\CX_a$ and $g$ copies of the sheaf $\Phi^{-1}(\ad)$ into these 
marked points. 

\subsection{} Note that an arbitrary $A$-curve may be degenerated into 
a curve considered in the previous example. Due to \ref{fus cor}, 
this determines the stalks of all our local systems (up to a non canonical  
isomorphism).

\section{Conformal blocks (higher genus)}

In this section we assume that $k=\BC$.  

\subsection{} Let us make the assumptions of \ref{charge}, \ref{level}. 
Consider the full subcategory $\MS\subset\CC$. Let us define the 
{\em regular object} $\bR_{\MS}$ by 
$$
\bR_{\MS}=
\oplus_{\mu\in\Delta_{\ell}}\ (\fL(\mu)\otimes\fL(\mu)^*)\in\MS^{\otimes 2}. 
$$
As in the previous section, we have a notion of $A$-collection 
$\{L_a;\bR_{\MS n}\}_A\ (A\in\Sew, \{L_a\}\in\MS^{A^0})$. The classical 
fusion structure on $\MS$, \cite{tuy}, defines for each $A$-collection 
as above, a local system 
$$
\langle\otimes_A\ \{L_a;\bR_{\MS n}\}\rangle_{\MS}
$$
on the moduli stack $\CM_{A;\lambda}$. 

We have $A=(\bA,N_A)$; let us define another object $A'=(\bA\cup\{*\},N_A)$. 
We have an obvious deleting $f_A:A'\lra A$. Given an $A$-collection as 
above, define an $A'$-collection in $\{L_a;L(2\rho_{\ell});\bR_n\}_{A'}$ in 
the category $\CC$. Using the equivalence $\Phi$, we transfer to $\CC$ the 
fusion structures defined in the previous section on $\FS$; 
we denote them $\langle\ \rangle^{(m)}_{\CC}$.  

The following theorem generalizes Theorem \ref{subq} to higher genus. 

\subsection{Theorem} 
\label{subquot} {\em For each $A\in\Sew$, the local system 
$f_A^*\langle\otimes_A\ \{L_a;{\em\bR}_{\MS n}\}\rangle_{\MS}$ on 
$\CM_{A';\lambda}$ is a canonical subquotient of the local system 
$\langle\otimes_{A'}\ \{L_a;L(2\rho_{\ell});{\em\bR}_n\}\rangle^{(0)}_{\CC}$.} 
$\Box$ 

\subsection{} Let us consider the special case $A=(A^0,\emp)$; 
for an $A$-collection $\{L_a\}_{A^0}$ in $\MS$, we have the classical  
local systems of conformal blocks $\langle\otimes_{A^0}\ \{L_a\}\rangle_{\MS}$ 
on $\CM_{A;\lambda}$. 

\subsection{Corollary} {\em The local systems 
$\langle\otimes_{A^0}\ \{L_a\}\rangle_{\MS}$ are semisimple. They carry a 
canonical non-degenerate Hermitian form.} 

In fact, the local system $f^*_A\langle\otimes_{A^0}\ \{L_a\}\rangle_{\MS}$ 
is semisimple by the previous theorem and by 
Beilinson-Bernstein-Deligne-Gabber, \cite{bbd} 6.2.5. The 
map of fundamental groupoids 
$f_{A*}:\tTeich_{A'}\lra\tTeich_A$ is surjective; therefore the initial 
local system is semisimple. 

The Hermitian form is defined in the same manner as in genus zero, cf.  
\ref{hermit}.

\end{document}